\documentclass{article}
\usepackage{amsmath}
\usepackage{amsfonts,amssymb}
\input amssym.def

\numberwithin{equation}{section}

\newcommand{\hA}{A_{\hbar}}
\newcommand{\dia}{\diamond}
\newcommand{\Der}{{\rm D e r}}
\newcommand{\pr}{{\rm p r}}
\newcommand{\intl}{{\rm i n t}}
\newcommand{\ext}{{\rm e x t}}
\newcommand{\Eu}{\, {\rm E u} \,}

\newcommand{\tnu}{\widetilde{\nu}}
\newcommand{\ttt}{\widetilde{t}}

\newcommand{\hpsi}{\hat{\psi}}
\newcommand{\tf}{\widetilde{f}}
\newcommand{\tg}{\widetilde{g}}
\newcommand{\tom}{\widetilde{\omega}}
\newcommand{\tcK}{\widetilde{\cal K}}
\newcommand{\tcS}{\widetilde{\cal S}}

\newcommand{\Ext}{{\rm E x t }}

\newcommand{\brarrow}{\succ\rightarrow}

\newcommand{\bbrarrow}{\succ\succ\rightarrow}
\newcommand{\bblarrow}{\leftarrow\prec\prec}

\newcommand{\SX}{{\cal S}X}

\newcommand{\Tp}{T_{poly}^{\bul}}

\newcommand{\cTp}{\cT^{\bul}_{poly}}

\newcommand{\cEb}{{\cal E}^{\bullet}}
\newcommand{\cAb}{{\cal A}^{\bullet}}
\newcommand{\Omb}{\Om^{\bul}}

\newcommand{\OmT}{\Om^{\bul}(\cT^{\bul}_{poly})}
\newcommand{\OmD}{\Om^{\bul}(\Cbu(\SX))}
\newcommand{\OmC}{\Om^{\bul}(\Cbd(\SX))}
\newcommand{\OmE}{\Om^{\bul}(\cE^{\bul})}

\newcommand{\Cbu}{C^{\bullet}}
\newcommand{\Cbd}{C_{\bullet}}

\newcommand{\Linf}{L_{\infty}}

\newcommand{\tV}{\widetilde{V}}

\newcommand{\la}{{\lambda}}
\newcommand{\h}{{\hbar}}
\newcommand{\bul}{{\bullet}}

\newcommand{\mA}{{\mathfrak{A}}}

\newcommand{\mb}{{\mathfrak{b}}}
\newcommand{\mg}{{\mathfrak{g}}}
\newcommand{\mh}{{\mathfrak{h}}}

\newcommand{\om}{{\omega}}

\newcommand{\Om}{{\Omega}}

\newcommand{\ga}{{\gamma}}

\newcommand{\G}{{\Gamma}}

\newcommand{\pa}{{\partial}}

\newcommand{\cK}{{\cal K}}

\newcommand{\cL}{{\cal L}}
\newcommand{\cD}{{\cal D}}
\newcommand{\cA}{{\cal A}}
\newcommand{\cG}{{\cal G}}

\newcommand{\cT}{{\cal T}}

\newcommand{\cS}{{\cal S}}
\newcommand{\cE}{{\cal E}}

\newcommand{\cO}{{\cal O}}

\newcommand{\bbC}{{\Bbb C}}
\newcommand{\bbR}{{\Bbb R}}
\newcommand{\bbZ}{{\Bbb Z}}

\newcommand{\La}{{\Lambda}}

\newcommand{\n}{{\nabla}}

\newcommand{\de}{{\delta}}

\date{}
\newtheorem{defi}{Definition}

\newtheorem{teo}{Theorem}
\newtheorem{cor}{Corollary}
\newtheorem{pred}{Proposition}

\author{Vasiliy Dolgushev}

\title{The Van den Bergh duality and
the modular symmetry of a Poisson variety}

\begin{document}

\maketitle

\begin{abstract}
We consider a smooth Poisson affine
variety with the trivial canonical bundle
over $\bbC$.
For such a variety the deformation quantization
algebra $A_{\h}$
enjoys the conditions of the Van den Bergh duality
theorem and the corresponding dualizing module
is determined by an outer
automorphism of $A_{\h}$ intrinsic to $A_{\h}$\,.
We show how this automorphism can be
expressed in terms of the modular class
of the corresponding Poisson variety.
We also prove that the
Van den Bergh dualizing module of the
deformation quantization algebra $\hA$ is free
if and only if the corresponding Poisson
structure is unimodular.
\end{abstract}

{ \large

\section{Introduction}
In seminal paper \cite{W} A. Weinstein, inspired
by Connes-Takesaki-Tomita theory of modular automorphisms
of a von Neumann algebra, introduced a notion of
modular class of a smooth real Poisson manifold $M$.
This class belongs to the first Poisson cohomology
group $HP^1(M)$, i.e. the space of Poisson vector fields
modulo Hamiltonian vector fields. Independently,
J.-L. Brylinski and G. Zuckerman
\cite{BG} introduced and studied the modular vector
field in the context of complex analytic
Poisson geometry. The modular class is a natural
generalization of the modular character of a Lie
algebra. This class also admits further generalizations
to Lie algebroids \cite{ELW}, to Lie-Rinehart
algebras \cite{Hueb}, and to $Q$-manifolds
\cite{Tomsk}.

In this paper we show that the deformation
quantization incarnation of the modular symmetry
of Poisson manifolds is
related to the Van den Bergh duality \cite{VB}
between Hochschild homology and Hochschild
cohomology. This relationship can be described 
in terms of Bursztyn-Waldmann bimodule 
quantization \cite{BW}. More precisely, a modular
vector field of a Poisson structure gives us a flat 
contravariant connection 
on the $A$-bimodule $A$, where $A$ is the algebra of 
functions on the Poisson variety. This contravariant 
connection can be quantized to a bimodule of the 
deformation quantization algebra $\hA$. The main result of 
the paper (see Theorem \ref{ONA}) states that the resulting  
bimodule is the Van den Bergh dualizing bimodule $\hA$\,.
In this paper we also show that
the Van den Bergh dualizing module of $A_\h$
is isomorphic to $A_{\h}$ if and only if the
corresponding Poisson structure $\pi$
is unimodular.

The organization of the paper is as follows.
In the second section we go over the notation
and recall some required results.
At the end of this
section we describe a construction which
produces out of a Poisson vector field
a derivation of the deformation quantization
algebra. The third section is devoted to the
definition of the modular class of a Poisson
structure in the algebraic setting.
In section 4 we define the modular
automorphism\footnote{In \cite{BZh}
K.A. Brown and J.J. Zhang call it the Nakayama
automorphism.} of a deformation
quantization algebra, formulate our
main result (Theorem \ref{ONA}), and prove
a useful technical proposition which we need 
in the proof of Theorem \ref{ONA}. 
Section 5 is devoted to
criterion of unimodularity and section 6
is devoted to the proof of Theorem \ref{ONA}.
In the concluding section we discuss some
results in literature related to Theorems \ref{ONA}
and \ref{unimod-teo}. In the Appendix we
discuss properties of Poisson, Hamiltonian
and log-Hamiltonian vector fields.

~\\
{\bf Acknowledgment.}
I would like to thank B. Enriquez, P. Etingof, V. Ginzburg,
Y. Kosmann-Schwarzbach, S. Launois and V. Rubtsov for discussions
of this topic. 
I am grateful to  B. Enriquez and V. Rubtsov for
their comments about various versions of the
manuscript. I am grateful to B. Enriquez for pointing 
out a significant gap in the original version of the paper. 
 The results of this paper were
presented in the seminars at several institutions.
These are the University of Angers, University of
Montpellier II, ETH in Z\"urich, and Washington University 
in St. Louis.
I would like to thank the participants of these
seminars for their questions and useful comments.
This work was in preparation when I was a visitor
at FIM at ETH in Z\"urich. I would like to thank
this institute for invitation and for perfect
working conditions. I also want to thank my landlord
in Z\"urich Kerry J. Hines-Randle for his hospitality.
I am partially supported by
the Grant for Support of Scientific
Schools NSh-8065.2006.2.

\section{Preliminaries}
For an associative algebra $B$
we denote by $B^{op}$ the algebra with the opposite
multiplication and by $B^e$ the enveloping algebra
$B^e =  B\otimes B^{op}$ of $B$\,. $\Der(B)$ denotes
the Lie algebra of derivations of $B$\,.
$\Cbd(B, N)$ is the Hochschild chain complex
of $B$ with coefficients in the
$B$-bimodule $N$
\begin{equation}
\label{chains}
\Cbd(B,N) = N \otimes B^{\otimes \, \bul}
\end{equation}
and $\Cbu(B,N)$ is the
Hochschild cochain complex of $B$
with coefficients in $N$
\begin{equation}
\label{cochains}
\Cbu(B, N) = H o m(B^{\otimes\, \bul} , N)\,.
\end{equation}
The Hochschild coboundary operator
is denoted by $\pa$ and the
Hochschild boundary operator is
denoted by $\mb$\,. Furthermore, we reserve
the notation $HH_{\bul}(B,N)$ for the
homology of the complex $(\Cbd(B,N), \mb)$ and
the notation $HH^{\bul}(B,N)$ for the cohomology
of the complex $(\Cbu(B,N),\pa)$

By convention, $\Cbd(B)$ is the Hochschild
chain complex $\Cbd(B,B)$ of $B$ with coefficients
in $B$:
\begin{equation}
\label{chains-B}
\Cbd(B) = C_{\bul}(B,B)\,.
\end{equation}
While $\Cbu(B)$ stands for the Hochschild
cochain complex of $B$ with coefficients $B$ and
with the shifted grading:
\begin{equation}
\label{rule}
\Cbu(B) = C^{\bul+1} (B,B)\,.
\end{equation}
In particular the lowest component of (\ref{rule})
is
$$
C^{-1}(B) = B\,.
$$

``DGLA'' always means a differential graded Lie algebra.
The arrow $\brarrow$ denotes an $\Linf$-morphism
of DGLAs,
the arrow $\bbrarrow$ denotes a
morphism of $\Linf$-modules,
and the notation
$$
\begin{array}{c}
L\\[0.3cm]
~~ \downarrow_{\,mod}\\[0.3cm]
M
\end{array}
$$
means that $M$ is a DGLA module
over the DGLA $L$\,. The symbol
$\circ$ always stands for the composition
of morphisms. $\h$ denotes the
deformation parameter.

Throughout this paper
$X$ is a smooth complex affine variety
of dimension $d$ with the trivial canonical
bundle. $A=\cO(X)$ denotes the algebra of regular
functions on $X$\,. $TX$ (resp. $T^*X$) denotes the
tangent (resp. cotangent) sheaf of $X$\,.
Since $X$ is smooth, irreducible components of
$X$ are exactly its connected components.
For this reason it is sufficient to formulate
all our proofs in the case of irreducible
variety.

We denote by $\G(\cG)$  the module
of global sections of the sheaf
$\cG$ and by  $\Omb(\cG)$ the module
of exterior forms with values in $\cG$\,.
In this paper we only consider quasi-coherent
sheaves of modules over $\cO_X$ and, since
$X$ is affine, we tacitly identify
all sheaves with the corresponding
modules of global sections.

It is well known that the
Hochschild cochain complex (\ref{rule})
carries the structure of a DGLA. The corresponding
Lie bracket (see Eq. (3.2) on page 45 in \cite{thesis})
was originally introduced by M. Gerstenhaber
in \cite{Ger}. We will denote this bracket by
$[,]_G$\,.

The Hochschild chain complex $\Cbd(B)$ (\ref{chains-B})
carries the structure\footnote{To be more precise,
it is the Hochschild chain complex with
the reversed grading $C_{-\bul}(B)$ which
carries this DGLA module structure.}
of a DGLA module over
the DGLA $\Cbu(B)$\,. We will denote the
action (see Eq. (3.5) on page 46 in \cite{thesis})
of cochains on chains by $R$\,.

Due to \cite{HKR}
the cohomology of the complex $\Cbu(A)$ of
the algebra $A$ of functions on $X$ is the
module $\Tp(X)$ of polyvector fields on $X$
with the shifted grading
\begin{equation}
\label{Tp}
\Tp(X) = \G(\wedge^{\bul+1}_{\cO_X} TX)\,,
\qquad
T^{-1}_{poly}(X) = \cO(X)\,.
\end{equation}
The homology of the complex $\Cbd(A)$
is the module of exterior forms
\begin{equation}
\label{cAb}
\cAb(X) = \G(\wedge^{\bul}_{\cO_X} T^*X)\,.
\end{equation}

$\Tp(X)$ is a graded Lie algebra
with respect the so-called Schouten-Nijenhuis
bracket $[,]_{SN}$
(see Equations (\ref{SN}), (\ref{SN1})) and
$\cAb(X)$ is a graded Lie algebra module
over $\Tp(X)$ with respect to the Lie derivative $\cL$
(\ref{Lie}). We will regard $\Tp(X)$
(resp. $\cAb(X)$) as the DGLA
(resp. the DGLA module)
with the vanishing differential.

We denote by $x^i$ local coordinates on $X$ and
by $y^i$ fiber coordinates in the tangent bundle
$TX$. Having these coordinates $y^i$ we
can introduce another local basis of exterior
forms $\{ dy^i \}$. We will use both bases
$ \{ dx^i \}$ and $\{ dy^i \}$. In particular,
the notation $\Omb(\cG)$
is  reserved for the module of
$d y$-exterior forms with values in
the sheaf $\cG$ while $\cAb(X)$ (\ref{cAb}) denotes the
module of $d x$-exterior forms.

$\SX$ is the formally completed symmetric algebra
of the cotangent bundle $T^*X$\,. Elements of
the algebra $\SX$ can be viewed as formal power series in
tangent coordinates $y^i$\,.
We regard $\SX$ as
the algebra over $A$\,. In particular,
$\Cbu(\SX)$ (resp. $\Cbd(\SX)$) is the Hochschild
cochain (resp. chain) complex
of $\SX$ over $A$.

As in \cite{thesis} the tensor product in
$$
C_k(\SX) =
\underbrace{\SX \hat{\otimes}_{A} \SX
\hat{\otimes}_{A}
\dots
\hat{\otimes}_{A}  \SX}_{k+1}
$$
is completed in the adic topology in fiber
coordinates $y^i$ on the tangent bundle $TX$\,.
Similarly, the Hochschild cochains of $\SX$
are substituted by the formal fiberwise
polydifferential operators (see Definition 12
on page 60 in \cite{thesis}).

The cohomology of the complex
$\Cbu(\SX)$ is the module $\cTp$ of fiberwise
polyvector fields (see Definition 11 on page 60
in \cite{thesis}).
The homology of the complex
$\Cbd(\SX)$ is the module $\cEb$ of fiberwise differential
forms (see page 62 in \cite{thesis}). These are
$d x$-forms with values in $\SX$\,.
As $\Tp(X)$, the module $\cTp$ is also equipped
with a Schouten-Nijenhuis bracket.
Similarly $\cEb$ forms a graded Lie
algebra module
over $\cTp$. We will use the same notation
for the bracket $[\,,\,]_{SN}$ and for the
Lie derivative $\cL$\,. However, we have to keep
in mind that, unlike on $\Tp(X)$ and $\cAb(X)$,
the operations $[\,,\,]_{SN}$ and $\cL$ on
$\cTp$ and $\cEb$ are $A$-linear.

In \cite{thesis} (see Theorem 4 on page 68) it
is shown (in the $C^{\infty}$ setting) that
the algebra $\Omb(\SX)$
can be equipped with a differential of the
following form
\begin{equation}
\label{DDD}
D = \n - \de + \La\,,
\end{equation}
where
\begin{equation}
\label{nabla}
\n = dy^i \frac{\pa}{\pa x^i} -
dy^i \G^k_{ij}(x) y^j \frac{\pa}{\pa y^k}\,,
\end{equation}
is a torsion free connection with
Christoffel symbols $\G^k_{ij}(x)$,
\begin{equation}
\label{delta}
\de = dy^i \frac{\pa}{\pa y^i}\,,
\end{equation}
and
$$
\La =\sum_{p=2}^{\infty}dy^k \La^j_{ki_1\dots i_p}(x)
y^{i_1} \dots
y^{i_p}\frac{\pa}{\pa y^j} \in \Om^1(\cT^0_{poly})\,.
$$
We refer to (\ref{DDD}) as the Fedosov differential.

In our setting it is also possible
to construct the Fedosov differential (\ref{DDD})
since $X$ is affine, and hence $TX$ can be equipped
with the desired connection.

Notice that $\de$ in (\ref{delta}) is also a differential on
$\Omb(\SX)$ and (\ref{DDD}) can be viewed as
deformation of $\de$ via the connection $\n$\,.

According to Proposition $10$ on
page $64$ in \cite{thesis} the modules
$\cTp$, $\Cbu(\SX)$,
$\cE^{\bul}$, and $\Cbd(\SX)$ are equipped with the
canonical action of the Lie algebra
$\cT_{poly}^0$ and this action is compatible
with the corresponding (DG) algebraic structures.
Using this action in
chapter $4$ of \cite{thesis} we extend
the Fedosov differential (\ref{DDD}) to
differentials on the DGLAs (resp. DGLA modules)
$\Omb(\cTp)$, $\Omb(\cE^{\bul})$,
$\Omb(\Cbd(\SX))$, and $\Omb(\Cbu(\SX))$\,.

Using acyclicity of the Fedosov differential (\ref{DDD})
in positive dimension
one constructs in \cite{thesis} embeddings
of DGLA
modules\footnote{See Eq. (5.1) on page 81
in \cite{thesis}.}
\begin{equation}
\begin{array}{ccc}
\Tp(X) &\stackrel{\la_T}{\,\longrightarrow\,} &
(\OmT, D, [,]_{SN})\\[0.3cm]
\downarrow_{\,mod}  & ~  &
\downarrow_{\,mod}\\[0.3cm]
\cAb(X)  &\stackrel{\la_{\cA}}{\,\longrightarrow\,} &
(\OmE, D),
\end{array}
\label{diag-T}
\end{equation}

\begin{equation}
\begin{array}{ccc}
(\OmD, D+\pa, [,]_{G}) &\stackrel{\,\la_D}{\,\longleftarrow\,}
& \Cbu(A)\\[0.3cm]
\downarrow_{\,mod}  & ~  &
\downarrow_{\,mod} \\[0.3cm]
(\OmC, D+\mb) &\stackrel{\,\la_C}{\,\longleftarrow\,} &
\Cbd(A),
\end{array}
\label{diag-D}
\end{equation}
and shows that these are quasi-isomorphisms of the
corresponding complexes.

Furthermore, using Kontsevich's and Shoikhet's
formality theorems for $\bbR^d$ \cite{K}, \cite{Sh}
in \cite{thesis} one constructs the following
commutative diagram
\begin{equation}
\label{diag-K-Sh}
\begin{array}{ccc}
(\OmT, D, [,]_{SN}) & \stackrel{\cK}{\brarrow} &
(\OmD, D+\pa, [,]_{G}) \\[0.3cm]
\downarrow_{\,mod}  & ~  &
\downarrow_{\,mod} \\[0.3cm]
(\OmE, D) & \stackrel{\cS}{\bblarrow} &
(\OmC, D+\mb)
\end{array}
\end{equation}
where $\cK$ is an $\Linf$ quasi-isomorphism of
DGLAs, $\cS$ is a quasi-isomorphism
of $\Linf$-modules over the
DGLA $(\OmT, D, [,]_{SN})$\,,
and the $\Linf$-module structure on $\OmC$ is obtained
by composing the quasi-isomorphism
$\cK$ with the DGLA modules structure $R$
(see Eq. (3.5) on p. 46 in \cite{thesis} for the
definition of $R$)\,.

~\\
{\bf Remark.} As in \cite{thesis} we use adapted versions
of Hochschild (co)chains for the algebra of functions $A$.
Thus, $\Cbu(A)$ is the complex of (algebraic)
polydifferential operators (see page 48 in \cite{thesis})
on $X$ and $\Cbd(A)$ is the sheaf of polyjets
as in Equation (3.15) on page 49 in \cite{thesis}.
In the algebraic setting these
complexes are quasi-isomorphic to the corresponding
genuine Hochschild complexes (\ref{chains-B}) and
(\ref{rule}) of $A$\,.

Let us recall the formulas of Cartan
calculus of exterior forms and polyvector
fields.

The Lie derivative
with respect to a polyvector field $\ga$ is
defined as the graded commutator
\begin{equation}
\label{Lie}
\cL_{\ga} = [d, i_{\ga}]
\end{equation}
of the de Rham differential $d$ and
the contraction $i_{\ga}$ with $\ga$\,.
The contraction is defined in the
obvious way for $\ga$ being a function or
a vector and then extended to an arbitrary
polyvector by the equation
\begin{equation}
\label{i}
i_{\ga_1} i_{\ga_2} = i_{\ga_1 \wedge \ga_2}\,,
\qquad \ga_1, \ga_2 \in \Tp(X)\,.
\end{equation}

The Schouten-Nijenhuis bracket $[\,,\,]_{SN}$
on $\Tp(X)$ is defined by the equations
\begin{equation}
\label{SN}
[a,b]_{SN} =0\,, \qquad [w, a]_{SN}= w(a)\,,
\qquad [w_1, w_2]_{SN} = [w_1,w_2]\,,
\end{equation}
\begin{equation}
\label{SN1}
[\ga_1,  \ga_2 \wedge  \ga_3 ]_{SN} =
[\ga_1, \ga_2]_{SN} \wedge \ga_3 +
(-1)^{|\ga_1|(|\ga_2|+1)}
\ga_2 \wedge [\ga_1, \ga_3]_{SN}\,,
\qquad \ga_i \in \Tp(X)\,.
\end{equation}
where $a, b$ are functions, $w, w_1, w_2$
are vector fields,
$[\,,\,]$ stands for the Lie bracket of vector
fields, and $|\ga_i|$ denotes the degree of the
polyvector $\ga_i$ in $\Tp(X)$ (\ref{Tp})\,.

The Lie derivative and the contraction
operation satisfy the following
equations
\begin{equation}
\label{Cartan}
[i_{\ga_1}, [d, i_{\ga_2}]] = i_{[\ga_1, \ga_2]_{SN}}\,,
\end{equation}
\begin{equation}
\label{Cartan1}
[\cL_{\ga_1}, \cL_{\ga_2} ]
 = \cL_{[\ga_1, \ga_2]_{SN}}\,.
\end{equation}

Let us also recall the Van den Bergh
duality theorem:
\begin{teo}[M. Van den Bergh, \cite{VB}]
\label{VdB}
If $B$ is a finitely generated bimodule 
coherent\footnote{An algebra is called bimodule coherent 
if every map between finite rank free $B$-bimodules 
has a finitely generated kernel (see Definition 3.5.1 
in \cite{Vitya}).}
algebra of finite Hochschild dimension $d$\,,
\begin{equation}
\label{Ext}
\Ext^m_{B^e} (B, B \otimes B) =
\begin{cases}
U_B  & {\rm if} \quad m = d\,,\\
0  & {\rm otherwise}\,,
\end{cases}
\end{equation}
and $U_B$ is an invertible\footnote{A $B$-bimodule $U$ is
called invertible if there is a $B$-bimodule $V$ such
that $U \otimes_B \, V \cong B$\,.}
$B$-bimodule then
for every $B$-bimodule $N$
$$
HH^{\bul}(B, N) \cong HH_{d - \bul}(B, U_B \otimes_B N)\,.
$$
\end{teo}
In Equation (\ref{Ext}) $B\otimes B$ is considered
as a bimodule over $B$ with respect
to the external $B$-bimodule structure.
It is the internal $B$-bimodule
structure which equips all the ${\rm E x t}$
groups $\Ext^{\bul}_{B^e} (B, B \otimes B)$
with a structure of $B$-bimodule.

We refer to $U_B$ as {\it the Van den Bergh dualizing
module} of $B$\,.

\subsection{Quantization of the derivations of
Poisson algebra}
\label{q-deriv}
Let $\pi_1$ be a Poisson structure on $X$ and
$\hA=(A[[\h]], *)$ be a deformation quantization
of $\pi_1$ in the sense of \cite{Bayen} and
\cite{Ber}. Let, furthermore,
\begin{equation}
\label{pi-h1}
\pi = \h \pi_1 + \h^2 \pi_2 + \dots
\in \h\, \G(\wedge^2 TX)[[\h]]
\end{equation}
be a representative of Kontsevich's class
of the deformation $\hA$\,.

Due to the Jacobi identity $[\pi, \pi]_{SN}=0$
we get a non-zero differential
$$
\pa_{\pi} = [\pi,\, ]_{SN}
$$
on $\Tp(X)[[\h]]$\,.
This differential was originally introduced
by Lichnerowicz in \cite{Lich} and
the cohomology of the complex
\begin{equation}
\label{Poisson-com}
(\Tp(X)[[\h]], [\pi,\, ]_{SN})
\end{equation}
is the desuspended
Poisson cohomology of $\pi$:
\begin{equation}
\label{HP}
HP^{\bul}(X,\pi) =
H^{\bul - 1} (\,\Tp(X)[[\h]], [\pi,\,]_{SN})\,.
\end{equation}
In particular, the zeroth Poisson cohomology
is exactly the $\bbC[[\h]]$-module of
Casimir functions of $\pi$ and the
first Poisson cohomology is the quotient
of Poisson vector fields of $\pi$ by
Hamiltonian vector fields.

In this subsection we construct
a $\bbC[[\h]]$-linear map
\begin{equation}
\label{cD}
w \to \cD_w : \G(TX)[[\h]]\cap \ker [\pi,\,]_{SN}
\to \Der (\hA)
\end{equation}
satisfying the following properties:
\begin{equation}
\label{cD-p0}
\cD_w = w \qquad mod \quad \h \,,
\end{equation}
\begin{equation}
\label{cD-p1}
[\cD_{w_1}, \cD_{w_2}] =
\cD_{[w_1, w_2]} + {\rm inner~ derivations}
\end{equation}
$$
\forall \quad w, w_1, w_2 \in
\G(TX)[[\h]]\cap \ker [\pi,\,]_{SN}\,.
$$

First, the Poisson structure $\pi$ (\ref{pi-h1})
lifts to a Maurer-Cartan element $\la_T(\pi)$
in $\OmT[[\h]]$ which is flat with respect to
the Fedosov differential $D$ (\ref{DDD})\,.
Using this element we extend the differential
$D$ on $\OmT[[\h]]$ to
\begin{equation}
\label{la-pi}
D + [\la_T(\pi),\,]_{SN} : \OmT[[\h]] \to (\OmT[[\h]])[1]\,,
\end{equation}
where $[1]$ denotes the shift of the total
degree by $1$\,.

Second, the star-product $*$ viewed as
an element in $C^1(A)[[\h]]$ lifts to
a $D$-flat cochain
in $C^1(\SX)[[\h]]$ and hence gives us a new
associative product on $\SX[[\h]]$
\begin{equation}
\label{diamond}
\dia = \la_D(*)
\end{equation}
compatible with the differential $D$ (\ref{DDD})\,.

Using this product we extend the original
differential $D + \pa$ on $\OmD[[\h]]$ to
\begin{equation}
\label{D-dia}
D + \pa_{\dia} : \OmD[[\h]]
\to (\OmD[[\h]])[1]\,,
\end{equation}
where $\pa_{\dia}$ is the Hochschild coboundary
operator corresponding to the new product
$\dia$ on $\SX[[\h]]$ and $[1]$ as above
denotes the shift of the total degree by $1$\,.

Next, following the
lines of section 5.3 in \cite{thesis} we
can construct the following chain of ($\Linf$)
quasi-isomorphisms of DGLAs
\begin{equation}
\label{chain}
(\,\Tp(X)[[\h]], [\pi,\,]_{SN})
\stackrel{\la_T}{\,\longrightarrow\,}
\OmT[[\h]] \stackrel{\tcK}{\,\brarrow\,}
\OmD[[\h]] \stackrel{\la_D}{\,\longleftarrow\,}
\Cbu(\hA)\,,
\end{equation}
where $\OmT[[\h]]$ (resp. $\OmD[[\h]]$)
carries the differential (\ref{la-pi})
(resp. the differential (\ref{D-dia})).

The maps $\la_T$ (\ref{diag-T})
and $\la_D$ (\ref{diag-D}) are
genuine morphisms of DGLAs and $\tcK$ is an
$\Linf$-quasi-isomorphism obtained from
$\cK$ in (\ref{diag-K-Sh}) by twisting
via the Maurer-Cartan element $\la_T(\pi)$
and adjusting by the action of the prounipotent
group corresponding to the Lie algebra
$$
\mg = \h \, C^0(\SX)[[\h]] \oplus
\h \, \Om^{1}(C^{-1}(\SX)) [[\h]]\,.
$$

Let $w\in \G(TX)[[\h]]$ be a Poisson vector
field of $\pi$\,. In other words, $w$ is a degree zero
cocycle in the complex $\Tp(X)[[\h]]$
with the differential $[\pi,\,]_{SN}$\,.
The map $\la_T$ in the chain (\ref{chain})
pulls $w$ to a degree zero cocycle $\la_T(w)$ in
the complex $\OmT[[\h]]$ with the differential
(\ref{la-pi}). Then, using the structure map
$\tcK_1$ of the first level of the $\Linf$-morphism
$\tcK$\,, we get a cocycle
\begin{equation}
\label{W}
W = \tcK_1 (\la_T(w))
\end{equation}
in the complex $\OmD[[\h]]$ with the
differential (\ref{D-dia})\,.

As a cocycle of degree zero, $W$ has two components:
\begin{equation}
\label{W01}
W = W_0 + W_1\,, \qquad W_0 \in
\Om^0(C^0(\SX))[[\h]]\,, \qquad W_1 \in \Om^1(\SX)[[\h]]\,.
\end{equation}
Recall that, because of the shift (\ref{rule}),
$C^0(\SX)= C^1(\SX,\SX)$ and $C^{-1}(\SX) =\SX$\,.

Since the Fedosov differential (\ref{DDD}) is
acyclic in positive exterior degree we can kill
the component $W_1$ by adding an exact term to $W$\,.
Namely, using the homotopy operator $\Phi$
(see Eq. (4.36) in chapter 4 of \cite{thesis})
for the Fedosov differential, we conclude that
the cocycle
\begin{equation}
\label{U0}
U = W_0 - \pa_{\dia} \Phi(W_1)
\in \Om^0(C^0(\SX))[[\h]]
\end{equation}
is cohomologous to $W$\,.
If we denote by $\pr_1$ the projection
$$
\pr_1 : \Om^0(C^0(\SX))[[\h]] \oplus
\Om^1(\SX)[[\h]]
 \to \Om^1(\SX)[[\h]]
$$
onto the degree $1$ exterior forms then
$U$ can be rewritten as
$$
U = W - (D+\pa_{\dia}) \Phi (\,\pr_1(W)\,)\,.
$$
Since $U$ has the only non-zero
component in the exterior degree $0$
the equation $(D+ \pa_{\dia})U=0$ is
equivalent to the pair of equations
$$
D U=0\,, \qquad \pa_{\dia} U = 0\,.
$$
On the other hand, by definition,
the map $\la_D$ in (\ref{chain}) identifies
$\Cbu(\hA)$ with $\Cbu(\SX) \cap \ker D$\,.
Hence, the equation
\begin{equation}
\label{cD-def}
\la_D(\cD_w ) = W - (D+\pa_{\dia})
\Phi(\, \pr_1(W)\, )
\end{equation}
defines an element $\cD_{w}$ in $C^1(\hA)$\,.
Furthermore, since $\pa_{\dia} U =0$ the
element $\cD_{w}$ is a derivation of $\hA$\,.
We take (\ref{W}) and (\ref{cD-def}) as
the equations defining the desired map $\cD$
(\ref{cD}).

Equation (\ref{cD-p0}) follows from the construction
while Equation (\ref{cD-p1}) follows from the
identity between structure maps $\tcK_1$ and
$\tcK_2$ of $\tcK$:
$$
[\tcK_1(\ga_1), \tcK_1(\ga_2)]_G -
\tcK_1([\ga_1, \ga_2]_{SN}) =
(D +\pa_{\dia}) \tcK_2(\ga_1, \ga_2) +
\tcK_2( D \ga_1, \ga_2) +
(-1)^{|\ga_1|} \tcK_2(\ga_1, D \ga_2)
+
$$
$$
\tcK_2([\la_T(\pi),\ga_1]_{SN}, \ga_2) +
 (-1)^{|\ga_1|} \tcK_2(\ga_1, [\la_T(\pi),\ga_2]_{SN})\,,
$$
$$
\ga_1, \ga_2\in \OmT[[\h]]\,,
$$
where $|\ga_1|$ is the degree of $\ga_1$\,.

\section{Classical modular symmetry}
Let $X$ be a smooth affine variety over $\bbC$
with a trivial canonical class and with
$\pi_1 \in \G(\wedge^2 TX)$ being
a Poisson structure.

Following M. Kontsevich \cite{K} equivalence classes of
star-products of $\pi_1$ on $X$ can be parameterized using
deformed Poisson structures. Those are formal power series
in $\h$ of sections of $\wedge^2 TX$
\begin{equation}
\label{pi-h}
\pi = \h \pi_1 + \h^2 \pi_2 + \dots
\end{equation}
starting with $\h \pi_1$ and
satisfying the Jacobi identity:
 \begin{equation}
\label{Jacobi}
[\pi, \pi]_{SN} =0\,.
\end{equation}

This motivates us to replace polyvector fields
and exterior forms by the corresponding
formal power series in $\h$. In particular,
it makes sense to consider the Poisson structure
(\ref{pi-h}), as well as its Poisson and
Hamiltonian vector fields from the very
beginning.

In this section we define
the modular class \cite{BG}, \cite{W}
of the Poisson structure
$\pi$ (\ref{pi-h}) on a smooth affine variety
$X$ with trivial canonical bundle.
Unlike in the real $C^{\infty}$ setting as
in \cite{W}, or complex analytic setting
as in \cite{BG}, in our case
it is no longer possible to take
$\log$ of nowhere vanishing functions.
This is why the modular class
is not a class in the first Poisson cohomology
of $\pi$\,. Instead it is defined by a Poisson
vector field considered modulo the so-called
log-Hamiltonian vector fields.
We introduce and discuss the latter fields in the
Appendix.

Let $\om$ be a nowhere vanishing regular
section of the canonical bundle $\wedge^{d}T^*X$
on $X$\,. Since the bundle $\wedge^{d}T^*X$
is trivial such a section always exists and
it is defined up to a multiplication by
a unit $f$ of $A= \cO_X$:
\begin{equation}
\label{om-f}
\om \,\, \to \,\, f \,\om\,.
\end{equation}

Let us define the
modular vector field $v\in \h \G(TX)[[\h]]$ of the
Poisson structure $\pi$ (\ref{pi-h})
by the equation:
\begin{equation}
\label{modular}
i_v \om = \cL_{\pi}\om\,.
\end{equation}
The vector field $v$ is well defined
since $\om$ is nowhere vanishing.

Let us show that $v$ is a Poisson
vector field of $\pi$ or in other words,
\begin{equation}
\label{pi-v}
[\pi, v]_{SN} = 0\,.
\end{equation}

Indeed, using (\ref{Cartan}),
(\ref{Cartan1}), (\ref{Jacobi}),
and (\ref{modular}) we get
$$
i_{[\pi, v]_{SN}} \om
= - i_{[v, \pi]_{SN}} \om =
 - i_v \cL_{\pi} \om - \cL_{\pi} i_v \om  =
 - i_v i_v \om -
\cL_{\pi} \cL_{\pi}\om  =0\,.
$$
Therefore, since $\om$ is nowhere
vanishing, the bivector field $[\pi,v]_{SN}$
should vanish.

If we alter the section $\om$ according
to (\ref{om-f}) the field $v$ differs as
\begin{equation}
\label{v-f}
v \,\, \to  \,\, v - f^{-1} [\pi, f]_{SN}\,.
\end{equation}
Thus, we see that $v$ is defined up to
adding the log-Hamiltonian vector fields
(See the Appendix.).

~\\
{\bf Remark.} We might as well use
the formal power series
of top degree forms in $\h$
\begin{equation}
\label{t-om}
\widetilde{\om} = \om_0 + \h\, \om_1 + \h^2 \, \om_2 +
\dots
\end{equation}
starting with a nowhere vanishing
form $\om_0$ instead of $\om$ in
the definition of the modular
vector field (\ref{modular}).

Following \cite{W} we say that
\begin{defi}
A Poisson structure $\pi$ (\ref{pi-h})
is called \emph{unimodular} if
its modular class is trivial.
\end{defi}
As we are working in the algebraic setting
the triviality of the modular class means
that a modular vector field of $\pi$ is a
log-Hamiltonian vector field of $\pi$.
In other words, there exists a unit $f$ in
$A[[\h]]$ such that
\begin{equation}
\label{f-om}
\cL_{\pi} \om  - f^{-1}i_{[\pi, f]_{SN}} \om =0\,,
\end{equation}
where $\om$ is the nowhere vanishing
top degree form which was used to
define the modular vector field.

Using the formulas of Cartan calculus (\ref{Lie}),
(\ref{Cartan}), and (\ref{Cartan1})
it is not hard to show that
Equation (\ref{f-om}) is equivalent to
$$
\cL_{\pi} (f \om) = 0\,.
$$
Thus, if $\pi$ is unimodular then
there exists a formal power series
$\tom$ (\ref{t-om}) of top degree forms
satisfying the
equation
\begin{equation}
\label{unimod11}
\cL_{\pi} \tom = 0\,.
\end{equation}
Adding to this observation
the obvious reverse statement we get
\begin{pred}
\label{unimod}
The Poisson structure $\pi$ (\ref{pi-h})
is unimodular if and only if
there exits a
formal power series of top degree forms
$$
\tom = \om_0 + \h\, \om_1 + \dots \in \Om^d(X)[[\h]]
$$
starting with
a nowhere vanishing form $\om_0$
such that
\begin{equation}
\label{unimod1}
\cL_{\pi} \tom = 0\,. \qquad \Box
\end{equation}
\end{pred}

~\\
{\bf Examples.} If the Poisson variety $(X,\pi_1)$
is symplectic then any Poisson structure $\pi$
of the form (\ref{pi-h}) is unimodular.
Indeed, in this case one can produce the
inverse $\eta$ for the series of bivectors
$$
\pi_1 + \h \,\pi_2 + \h^2 \, \pi_3 + \dots\,.
$$
Due to (\ref{Jacobi}) $\eta$ is a formal
series in $\h$ of closed two forms starting
with $\pi_1^{-1}$ and
$\wedge^{d/2} \eta$ is
the desired volume form (\ref{t-om})
satisfying Equation (\ref{unimod1}).
An example of a Poisson structure with
a non-trivial modular class
is provided by the
dual space $\mh^*$ of a finite dimensional
Lie algebra $\mh$ with a non-zero modular
character.

\section{Quantum modular symmetry}
Let $\hA=(A[[\h]], *)$ be a deformation
quantization algebra of the Poisson
variety $(X,\pi_1)$ in the sense of
\cite{Bayen} and \cite{Ber}. We regard
$\hA$ as an algebra over the ring $\bbC[[\h]]$\,.

\begin{pred}
\label{VdB-OK}
The algebra $\hA$ has the Hochschild dimension 
$d = d i m\, X$\,. It satisfies the conditions 
of Theorem \ref{VdB}. The dualizing bimodule of $\hA$
is isomorphic to 
$\hA \, \nu$, where $\nu$ is an automorphism 
of $\hA$ satisfying the property
$$
\nu = Id \qquad m o d \quad \h\,.
$$
\end{pred}
{\bf Proof.} Using the standard arguments of deformation theory
(see Proposition 6 in \cite{Pasha} or step 2
in the proof of Theorem 4.1 in \cite{Pasha-O})
one can show
that the algebra $\hA$ (over $\bbC[[\h]]$)
has the same Hochschild dimension as the algebra $A$,
\begin{equation}
\label{k-neravno-d}
\Ext^k_{\hA^e}
(\hA, \hA \otimes \hA ) = 0
\end{equation}
for every $k\neq d$, and there
is an automorphism $\nu$ of $\hA$
such that
$$
\nu = Id \qquad m o d \quad \h\,,
$$
and
\begin{equation}
\label{Ext-d}
 \Ext^d_{\hA^e}
(\hA, \hA\otimes \hA ) \cong
\hA \, \nu
\end{equation}
as bimodules over $\hA$\,.
Here $\hA\otimes \hA$ is considered
as a bimodule over $\hA$ with respect
to the external $\hA$-bimodule structure.
It is the internal $\hA$-bimodule
structure which equips all the ${\rm E x t}$
groups $\Ext^{\bul}_{\hA^e} (\hA, \hA\otimes \hA)$
with a structure of $\hA$-bimodule.

$\hA$ is finitely generated since so is $A$\,.
Furthermore, 
it is not hard to show that bimodule coherence 
is stable under deformations we consider and the 
bimodule coherence of the commutative algebra 
$A$ follows from the fact that $A\otimes_{\bbC} A$
is Noetherian.  

The proposition is proved. $\Box$

Notice that, $\nu$ is defined up to
a composition with an inner automorphism.

\begin{defi}
We call the outer automorphism
$\nu$ \emph{the modular automorphism} of $\hA$\,.
\end{defi}

This definition is motivated
by the main result of this
paper which relates the outer
automorphism $\nu$ of the deformation
quantization algebra $\hA$ to
the modular class of the Poisson
structure $\pi$ (\ref{pi-h1}).
This result can be formulated as
follows:
\begin{teo}
\label{ONA}
Let $\hA$ be a deformation quantization
algebra of $(X, \pi_1)$ and let $\pi$
(\ref{pi-h1}) be a representative of Kontsevich's
class of $A_\h$\,.
If $v$ is a modular vector field of
the Poisson structure $\pi$ and
$\cD_v$ is a derivation of $A_\h$
constructed from $v$ via Kontsevich's
formality  theorem then the modular
automorphism $\nu$ of $A_\h$ is equal
to $\exp(\cD_v)$ up to an inner automorphism.
\end{teo}
{\bf Remark.}
It is Corollary \ref{sense} at the end
of the Appendix which implies that the
automorphism $\exp(\cD_v)$ in the above theorem
is well defined up to an inner automorphism.

In paper \cite{BW} H. Bursztyn and 
S. Waldmann showed that the semiclassical limit of 
an invertible bimodule over a deformation quantization 
algebra is a flat contravariant connection on an invertible 
bimodule of the corresponding Poisson algebra. 
They proved that, under certain conditions, such a 
contravariant connection on an invertible bimodule over $A$
can be quantized to an invertible bimodule over $\hA$\,. 
On the other hand every Poisson vector field 
defines a flat contravariant connection of the $A$-bimodule
$A$ and two such connections are isomorphic if and 
only if the corresponding Poisson vector fields 
differ by a log-Hamiltonian vector field
(see Lemma 4.8 in \cite{BW}). In terms of \cite{BW}
Theorem \ref{ONA} can formulated\footnote{The authors 
of \cite{BW} work in the category of $C^{\infty}$ real manifolds. 
For this reason we decided to formulate Theorem 
\ref{ONA} without referring to the quantization of 
H. Bursztyn and S. Waldmann \cite{BW}.} 
as follows: 
the deformation quantization of the contravariant 
connection corresponding to the modular vector 
field of a Poisson structure is the Van den Bergh 
dualizing bimodule of the deformation quantization 
algebra.

The proof of Theorem \ref{ONA} is given in
section \ref{Proof}\,. In the remainder of
this section we would like to prove an important
technical result which we will need
later on.

Let $w\in \h \, \G(TX)[[\h]]$ be a Poisson
vector field of $\pi$ (\ref{pi-h1}). Using
the map $\cD$ (\ref{cD}) we produce a derivation $\cD_{w}$
of the algebra $\hA$\,. Since the map $\cD$ is
$\bbC[[\h]]$-linear $\cD_w = 0 \, mod \, \h$\,.
Hence, it makes sense to speak about the
following automorphism
\begin{equation}
\label{phi-w}
\phi_w = \exp(\cD_w)
\end{equation}
of $\hA$\,.

Using this automorphism we define an
associative product $\cdot$
on $\hA[t,t^{-1}]$ by the following
rule:
\begin{equation}
\label{w}
(a t^n ) \cdot (b t^m) =
a * \phi_w^n (b) t^{n+m}\,,
\end{equation}
where $a,b \in \hA$ and $n,m$ are
arbitrary integers.

This product can be easily extended to 
a star-product on $A[t,t^{-1}][[\h]]$\,. 
In this way we obtain a deformation quantization
of the affine variety $X\times \bbC^{\times}$
with the following Poisson structure:
\begin{equation}
\label{pi-w-1}
\pi^w_1 = \pi_1 + t\pa_{t} \wedge w_1\,,
\end{equation}
where $\pi_1$ (resp. $w_1$) is the first coefficient
in the expansion of $\pi$ (\ref{pi-h1})
(resp. expansion of $w$) in $\h$ and
$\bbC^{\times} = \bbC\setminus \{ 0\}$\,.

The algebra $\hA[t,t^{-1}]$ with the product 
$\cdot$ (\ref{w}) is a proper subalgebra
of the deformation quantization algebra 
$(A[t,t^{-1}][[\h]], \cdot)$ of $X\times \bbC^{\times}$\,.
Thus, in order to apply the results of deformation 
quantization to $\hA[t,t^{-1}]$
we need to impose certain restrictions on
the formal power series of polyvector fields and 
exterior forms on $X\times \bbC^{\times}$\,.

We will denote by $\Tp(X\times \bbC^{\times})[[\h]]_{restr}$
the graded Lie algebra of formal power series of 
polyvector fields on $X\times \bbC^{\times}$ whose componets 
have bounded powers in $t$\,. Similarly, 
$\cAb(X\times \bbC^{\times})[[\h]]_{restr}$ will denote the
graded $\Tp(X\times \bbC^{\times})[[\h]]_{restr}$-module of 
formal power series of exterior forms on $X\times \bbC^{\times}$
whose components have bounded powers in $t$\,.

The bivector field 
\begin{equation}
\label{pi-w}
\pi^w  = \pi +  t\pa_{t} \wedge w
\end{equation}
obviously belongs to $\Tp(X\times \bbC^{\times})[[\h]]_{restr}$\,.
Furthermore, since $w$ is a Poisson vector field for $\pi$\,, the 
bivector
(\ref{pi-w}) satisfies the Jacobi identity: 
$$
[\pi^w , \pi^w ]_{SN} = 0\,.
$$ 
Hence $\pi^{w}$ equips the graded Lie algebra $\Tp(X\times \bbC^{\times})[[\h]]_{restr}$ and its module 
$\cAb(X\times \bbC^{\times})[[\h]]_{restr}$ with differentials 
$[\pi^w , \,\, ]$ and $\cL_{\pi^w}$\,, respectively.

We claim that
\begin{pred}
\label{dlia-hA-t}
There exists a chain of ($\Linf$-) quasi-isomorphisms 
connecting the DGLA module 
$(\Tp(X\times \bbC^{\times})[[\h]]_{restr}\,,\, \cAb(X\times \bbC^{\times})[[\h]]_{restr})$ with the differentials 
$[\pi^w , \,\, ]$ and $\cL_{\pi^w}$
to the DGLA module 
$(\Cbu(\hA[t,t^{-1}]), \Cbd(\hA[t,t^{-1}]))$ of 
the Hochschild (co)chains of the algebra $\hA[t,t^{-1}]$
with the product (\ref{w})
\end{pred}
{\bf Proof.} The idea of the proof is to start 
with the corresponding chain of 
($\Linf$-) quasi-isomorphisms for the DGLA module 
$(\Cbu(A[t,t^{-1}][[\h]]), \Cbd(A[t,t^{-1}][[\h]]))$
of Hochschild (co)\-chains of the deformation quantization 
algebra $A[t,t^{-1}][[\h]]$ and show that this chain 
can be restricted to a chain between desired DGLA 
modules. In doing this, a very important role will 
be played by the Euler field 
\begin{equation}
\label{Euler}
\Eu = t \frac{\pa}{\pa t}\,.
\end{equation}
In particular, we will use the following 
Fedosov differential  
on $\Omb(\cS(X\times \bbC^{\times}))$:
\begin{equation}
\label{DDD-t}
D^t = d y_t (t \frac{\pa}{\pa t} - \frac{\pa}{\pa y_t})
+ D\,,
\end{equation}
where $y_t$ is an auxiliary formal variable and 
$D$ is the Fedosov differential (\ref{DDD})
on $X$\,.

Replacing $X$ by $X\times \bbC^{\times}$
and $D$ by $D^t$ in
diagrams (\ref{diag-T}), (\ref{diag-D}), and
(\ref{diag-K-Sh}) we get the chain of
($\Linf$-) quasi-isomorphisms of DGLA modules: 
$$
\begin{array}{cccc}
\Tp(X \times \bbC^{\times}) &\stackrel{\la^w_T}{\,\longrightarrow\,} &
(\Omb(X \times \bbC^{\times}, \cT^{\bul}_{poly}), D^t, [,]_{SN})
&  \stackrel{\cK}{\brarrow}  \\[0.3cm]
\downarrow_{\,mod}  & ~  &
\downarrow_{\,mod} & ~ \\[0.3cm]
\cAb(X \times \bbC^{\times})  &\stackrel{\la^w_{\cA}}{\,\longrightarrow\,} &
(\Omb(X \times \bbC^{\times}, \cE^{\bul}), D^t) 
&  \stackrel{\cS}{\bblarrow} ,
\end{array}
$$
\begin{equation}
\begin{array}{cccc}
 \stackrel{\cK}{\brarrow} &
(\Omb(\Cbu(\cS(X\times \bbC^{\times}))), D^t + \pa, [,]_{G}) &\stackrel{\,\la^w_D}{\,\longleftarrow\,}
& \Cbu(A[t,t^{-1}])\\[0.3cm]
~ & \downarrow_{\,mod}  & ~  &
\downarrow_{\,mod} \\[0.3cm]
 \stackrel{\cS}{\bblarrow} &
(\Omb(\Cbd(\cS(X\times \bbC^{\times}))), D^t + \mb) &\stackrel{\,\la^w_C}{\,\longleftarrow\,} &
\Cbd(A[t,t^{-1}])\,,
\end{array}
\label{t-chain}
\end{equation}
where $\la^w_{T}$, $\la^w_{\cA}$, $\la^w_{D}$, and 
$\la^w_{C}$ are the corresponding versions of 
$\la_{T}$, $\la_{\cA}$, $\la_{D}$, and 
$\la_{C}$ in (\ref{diag-T}) and (\ref{diag-D}) 
for $X\times \bbC^{\times}$\,.

Twisting the chain (\ref{t-chain}) by the Maurer-Cartan 
element $\pi^w$ (\ref{pi-w}) we get the following 
chain of ($\Linf$-) quasi-isomorphisms of DGLA modules: 
$$
\begin{array}{cccc}
\Tp(X \times \bbC^{\times})[[\h]] &
\stackrel{\la^w_T}{\,\longrightarrow\,} &
(\Omb(X \times \bbC^{\times}, \cT^{\bul}_{poly})[[\h]], 
D^t +[\la^w_T(\pi^w),\,], [,]_{SN})
&  \stackrel{\tcK}{\brarrow}  \\[0.3cm]
\downarrow_{\,mod}  & ~  &
\downarrow_{\,mod} & ~ \\[0.3cm]
\cAb(X \times \bbC^{\times})[[\h]]  &\stackrel{\la^w_{\cA}}{\,\longrightarrow\,} &
(\Omb(X \times \bbC^{\times}, \cE^{\bul})[[\h]], 
D^t + \cL_{\la^w_T(\pi^w)}) 
&  \stackrel{\tcS}{\bblarrow} ,
\end{array}
$$
\begin{equation}
\begin{array}{cccc}
 \stackrel{\tcK}{\brarrow} &
(\Omb(\Cbu(\cS(X\times \bbC^{\times}))), D^t + \pa_{\dia_w}, [,]_{G}) &\stackrel{\,\la^w_D}{\,\longleftarrow\,}
& \Cbu(A[t,t^{-1}][[\h]], *_w)\\[0.3cm]
~ & \downarrow_{\,mod}  & ~  &
\downarrow_{\,mod} \\[0.3cm]
 \stackrel{\tcS}{\bblarrow} &
(\Omb(\Cbd(\cS(X\times \bbC^{\times}))), D^t + \mb_{\dia_w}) &\stackrel{\,\la^w_C}{\,\longleftarrow\,} &
\Cbd(A[t,t^{-1}][[\h]], *_w)\,,
\end{array}
\label{t-chain-h}
\end{equation}
where $\Tp(X \times \bbC^{\times})[[\h]]$ and 
$\cAb(X \times \bbC^{\times})[[\h]]$ are considered with 
the differentials $[\pi^w,\,\,]_{SN}$ and $\cL_{\pi^w}$
respectively, $*_w$ is a star-product on $A[t,t^{-1}][[\h]]$\,,
and $\dia_w$ is its lift to the algebra 
$\cS(X\times \bbC^{\times})[[\h]]$: 
\begin{equation}
\label{dia-w}
\dia_w = \la^w_D(*_w)\,.
\end{equation}

Let us note that the lift of the Euler field $\Eu$ 
(\ref{Euler}) to $\Om^0(X\times \bbC^{\times}, \cT^0_{poly}) 
\cap \ker D^t$ is given by 
\begin{equation}
\label{Eu-lift}
\la^w_T(\Eu) = \pa_{y_t}\,.
\end{equation}

It is not hard to see that the maps in the 
chain (\ref{t-chain}) are compatible with 
the Euler field 
$\Eu$ (\ref{Euler}) in the following sense
$$
\la^w_T([\Eu, \ga]_{SN}) = [\pa_{y_t}, \la^w_T(\ga)]_{SN}\,, 
\qquad 
\la^w_{\cA}(\cL_{\Eu}\, \eta ) = \cL_{\pa_{y_t}}\, \la^w_{\cA}(\eta ) 
$$
\begin{equation}
\label{Eu-compat}
\begin{array}{c}
\displaystyle
[\pa_{y_t}, \cK_n(\ga_1, \dots, \ga_n)]_{G} = 
\sum_{j=1}^{n} \cK_n(\ga_1, \dots, \ga_{j-1}, 
 [\pa_{y_t}, \ga_j]_{SN}\,,\ga_{j+1}, \dots,  \ga_n) \\[0.3cm]
\displaystyle
\cL_{\pa_{y_t}} \cS_n(\ga_1, \dots, \ga_n ; c) 
= \sum_{j=1}^{n} \cS_n(\ga_1, \dots, \ga_{j-1}, 
 [\pa_{y_t}, \ga_j]_{SN}\,,\ga_{j+1}, \dots,  \ga_n; c) + 
\end{array}
\end{equation}
$$
 +  \cS_n(\ga_1, \dots, \ga_n ; R_{\pa_{y_t}}\, c) 
$$
$$
\la^w_D([\Eu, P]_{G}) = [\pa_{y_t}, \la^w_D(P)]_{G}\,, 
\qquad 
\la^w_{C}(R_{\Eu}\, b ) = R_{\pa_{y_t}}\, \la^w_{C}(b)\,, 
$$
where $\ga\in \Tp(X\times \bbC^{\times})$\,,  
$\eta \in \cAb(X\times \bbC^{\times})$\,, 
$\ga_j \in \Omb(X \times \bbC^{\times}, \cT^{\bul}_{poly})$\,, 
$c\in \Omb(\Cbd(\cS(X\times \bbC^{\times})))$\,, 
$P \in \Cbu(A[t,t^{-1}])$\,, and $b \in \Cbd(A[t,t^{-1}])$\,.

On the other hand $[\Eu, \pi^w]_{SN} = 0$. Hence the
maps in the chain (\ref{t-chain-h}) are also compatible 
with the Euler field $\Eu$ and $*_w$ can be chosen 
to satisfy the property
$$
\Eu(a_1 *_w a_2) = \Eu (a_1)\, *_w\, a_2 + 
a_1\, *_w \Eu(a_2)\,, \qquad 
a_1, a_2 \in A[t,t^{-1}][[\h]]\,.
$$

Thus, since the Euler field ``counts the powers'' in $t$\,, 
the chain (\ref{t-chain-h}) restricts to the chain of 
($\Linf$-) quasi-isomorphisms which connects the DGLA module 
$(\Tp(X\times \bbC^{\times})[[\h]]_{restr}\,,\, \cAb(X\times \bbC^{\times})[[\h]]_{restr})$ with the differentials 
$[\pi^w , \,\, ]$ and $\cL_{\pi^w}$
to the DGLA module 
$(\Cbu(A[[\h]][t,t^{-1}])$\,, $ \Cbd(A[[\h]][t,t^{-1}]))$ of 
the Hochschild (co)chains of the algebra $A[[\h]][t,t^{-1}]$
with the product $*_w$\,. 

Our purpose now is to
prove that the algebra $\hA[t,t^{-1}]$
with the product (\ref{w})
is isomorphic to the algebra $A[[\h]][t,t^{-1}]$
with the product $*_w$\,.

It is obvious that
\begin{equation}
\label{la-T-w-la-T}
\la^w_T \Big|_{A} = \la_T\,.
\end{equation}
Furthermore, $w$, $t$, $\pi^w$
give us the following $D^t$-flat sections
of $\cTp[[\h]]$ on $X\times \bbC^{\times}$
\begin{equation}
\label{sections}
\la^w_T(w) = \la_T(w)\,, \qquad
\la^w_T(t) = t\, e^{y_t}\,,
\end{equation}
and
\begin{equation}
\label{lift-pi-w}
\la^w_T(\pi^w) = \la_T(\pi) + \pa_{y_t} \wedge \la_T(w)
\end{equation}
where, by abuse of notation,
we denote by $\la_T(w)$ (resp. $\la_T(\pi)$) the
natural lift of the fiberwise vector
(resp. bivector) on $TX$ to $T (X\times \bbC^{\times})$

Notice that
the vector field $w$ (on $X\times \bbC^{\times}$)
is log-Hamiltonian
with respect to $\pi^w$ (\ref{pi-w})\,.
More precisely,
\begin{equation}
\label{w-log}
w = - t^{-1} [\pi^w, t]_{SN}\,.
\end{equation}
Therefore, due to Theorem \ref{lift-log-Ham}
there exists a function
$\ttt\in A[t, t^{-1}][[\h]]$ such that
\begin{equation}
\label{ttt-t}
\ttt = t \qquad mod \quad \h\,,
\end{equation}
and for every $a\in A[t, t^{-1}][[\h]]$
\begin{equation}
\label{ttt-w}
\exp(\cD^t_w)\, (a)  = \ttt *_w\, a *_w \,\ttt^{-1}\,,
\end{equation}
where the inverse $\ttt^{-1}$ is taken in the algebra
$(A[t, t^{-1}][[\h]], *_w)$
and $\cD^t_w$ denotes the derivation of
this algebra corresponding to
the Poisson vector field $w$ on
$X\times \bbC^{\times}$\,.

The compatibility of the chain of maps (\ref{t-chain-h})
with the Euler field $\Eu$ (\ref{Euler}) implies that 
$\ttt$ as well as its inverse $\ttt^{-1}$ belongs to 
the subalgebra $A[[\h]][t, t^{-1}]\subset A[t, t^{-1}][[\h]]$\,.

The desired isomorphism $\tau$ from $\hA[t, t^{-1}]$
to $(A[[\h]][t, t^{-1}], *_w)$ is then
defined on generators as follows
\begin{equation}
\label{isom}
\tau(a) = a\,, \qquad \tau(t) = \ttt\,,
\end{equation}
where $a\in A$\,.

It is Equation (\ref{ttt-t}) which implies
that (\ref{isom}) is indeed an isomorphism
of $\bbC[[\h]]$-modules.

To show that (\ref{isom}) is a map of algebras
we need to check that
\begin{equation}
\label{check1}
a *_w \, b = a * b\,, \qquad
\forall \quad  a,b \in A[[\h]]\,,
\end{equation}
and
\begin{equation}
\label{check11}
\ttt *_w\, a *_w \,\ttt^{-1} =
\exp(\cD_w)\, (a)\,,
\qquad  \forall \quad  a \in A[[\h]]\,.
\end{equation}

Equation (\ref{check1}) follows from
(\ref{la-T-w-la-T}),
(\ref{lift-pi-w}) and the identity
$$
\cK_n(\la_T(\pi) + \pa_{y_t} \wedge \la_T(w),
\la_T(\pi) + \pa_{y_t} \wedge \la_T(w), \dots,
\la_T(\pi) + \pa_{y_t} \wedge \la_T(w) )(s_1,s_2) =
$$
$$
\cK_n(\la_T(\pi),
\la_T(\pi), \dots,
\la_T(\pi))(s_1,s_2)\,,
$$
where $\cK_n$ are the structure maps of the
quasi-isomorphism $\cK$ in (\ref{t-chain})
and $s_1, s_2$ are sections of
$\cS(X\times \bbC^{\times})$ satisfying
the equations
$$
\frac{\pa}{\pa y_t} s_1 =
\frac{\pa}{\pa y_t} s_2 = 0\,.
$$

To prove (\ref{check11}) it suffices to
show that
\begin{equation}
\label{cD-cD-t}
\cD^t_w (a) = \cD_w(a)\,, \qquad
\forall \qquad a\in A[[\h]]\,,
\end{equation}
where in the left hand side
$w$ is viewed as a vector field
on $X\times \bbC^{\times}$ and
in the right hand side $w$ is viewed
as a vector field on $X$\,.

Again (\ref{cD-cD-t}) follows from
(\ref{la-T-w-la-T}), (\ref{sections}),
(\ref{lift-pi-w})
and the identity
$$
\cK_{n+1}(\la_T(\pi) + \pa_{y_t} \wedge \la_T(w),
\la_T(\pi) + \pa_{y_t} \wedge \la_T(w), \dots,
\la_T(\pi) + \pa_{y_t} \wedge \la_T(w), \la_T(w))(s) =
$$
$$
\cK_{n+1}(\la_T(\pi), \la_T(\pi), \dots,
\la_T(\pi), \la_T(w))(s)\,,
$$
$\cK_n$ are, as above, the structure maps of the
quasi-isomorphism $\cK$ in (\ref{t-chain})
 and $s$ is a section of
$\cS(X\times \bbC^{\times})$ satisfying
the equation
$$
\frac{\pa}{\pa y_t} s = 0\,.
$$

This concludes the proof of Proposition \ref{dlia-hA-t}. $\Box$

\section{Criterion of unimodularity}
Let, as above, $X$ be a smooth complex
affine variety with trivial canonical bundle.
$d i m\, X = d$, $\pi_1$ is a Poisson structure
on $X$, and $\hA$ is a deformation quantization
algebra of $\pi_1$ with Kontsevich's class
represented by the formal series
$\pi$ (\ref{pi-h1}).

We have
\begin{teo}
\label{unimod-teo}
The Van den Bergh dualizing module
\begin{equation}
\label{U}
U_{\hA} = \Ext^d_{\hA^e}
(\hA, \hA \otimes \hA )
\end{equation}
of $\hA$ is isomorphic to $\hA$ as a bimodule
if and only if the Poisson structure
$\pi$ (\ref{pi-h1}) is unimodular.
\end{teo}
Although the ``if'' part of this statement is
an immediate corollary of Theorem \ref{ONA}\,,
we give here an independent proof of both
implications because we will use
Theorem \ref{unimod-teo} in full
generality in the proof of Theorem \ref{ONA}.

~\\
{\bf Proof of Theorem
\ref{unimod-teo}.} Let us prove the implication
$\Rightarrow$\,.

Due to \cite{VB} the condition
$U_{\hA} \cong \hA$ implies that
we have an isomorphism
\begin{equation}
\label{iso}
V : HH^0(\hA, \hA) \stackrel{\sim}{\to} HH_{d}(\hA, \hA)\,.
\end{equation}
On the other hand, the formality theorems
for Hochschild (co)chains \cite{FTHC}, \cite{K}, \cite{Sh}
provide us with isomorphisms
\begin{equation}
\label{mu-co}
\mu_1 : HP^0(X,\pi)  \stackrel{\sim}{\to} HH^0(\hA,\hA)\,,
\end{equation}
\begin{equation}
\label{mu-ch}
\mu_2 : HH_{d}(\hA, \hA)
\stackrel{\sim}{\to} HP_{d}(X, \pi)\,,
\end{equation}
where $HP^{\bul}(X, \pi)$ denotes the Poisson
cohomology (\ref{HP}),
and $HP_{\bul}(X,\pi)$ denotes the
Poisson homology \cite{JLB}, \cite{Koszul}\,.
The latter is, by definition, the homology of the complex
$\cAb(X)[[\h]]$ of exterior forms on $X$
with the differential $\cL_{\pi}$:
\begin{equation}
\label{HPd}
HP_{\bul}(X,\pi) =
H_{\bul} (\cAb(X)[[\h]], \cL_{\pi})\,.
\end{equation}

Taking the class $[1]\in HP^0(X,\pi)$ represented by
the function $1\in A= \cO(X)$ we get the class
$$
\mu_2 \circ V \circ \mu_1([1]) \in HP_{d}(X,\pi)\,.
$$
This class is represented by a formal power
series:
\begin{equation}
\label{tom}
\tom =  \om_0 + \h \om_1 + \dots \in
\G(\wedge^d T^* X)[[\h]]
\end{equation}
such that
\begin{equation}
\label{pi-tom}
\cL_{\pi} \om = 0\,.
\end{equation}

Thus, in view of Proposition \ref{unimod},
it suffices to prove that $\om_0$ is
nowhere vanishing. For this we consider
the Van den Bergh isomorphism (\ref{iso})
modulo $\h$\,.

Indeed, modulo $\h$ the map (\ref{iso})
gives us the Van den Bergh isomorphism
\begin{equation}
\label{iso0}
V_0 : HH^0(A,A) \stackrel{\sim}{\to} HH_{d}(A,A)
\end{equation}
for the commutative algebra $A = \cO(X)$.

Due to Hochschild-Kostant-Rosenberg theorem
\cite{HKR} the map $V_0$ is an isomorphism
from $A$ onto $A$\,.
This observation immediately implies that
$\om_0$ is a nowhere vanishing form and the
implication $\Rightarrow$ follows.

Let us prove the implication $\Leftarrow$.
Since $\pi$ is unimodular Proposition \ref{unimod}
implies that there exists a formal power
series (\ref{tom}) of top degree forms
satisfying (\ref{pi-tom}) and
starting with a nowhere vanishing form
$\om_0$\,.

Isomorphism (\ref{mu-ch}) provides us
with a class $\mu_2^{-1}([\tom])$ in
$HH_{d}(\hA,\hA)$\,, where $[\tom]$ is the
class of the cycle $\tom$ in the
Poisson chain complex.

On the other hand the Van den Bergh theorem
gives us an isomorphism
\begin{equation}
\label{iso-nu}
\tV : HH_d(\hA,\hA)
\stackrel{\sim}{\to} HH^{0}(\hA, \nu^{-1}\, \hA )\,.
\end{equation}
where $\nu$ is the modular automorphism
of $\hA$\,.

Let us denote by $b$ a representative
of the image
$\tV(\mu^{-1}_2([\tom]))$ of the
class $\mu^{-1}_2([\tom])$ under the map
$\tV$\,. By definition, $b\in \hA$, and hence
is a formal power series
$$
b = b_0 + \h b_1 + \h^2 b_2 \dots\,.
$$

Considering the isomorphisms $\mu_2$ (\ref{mu-ch})
and $\tV$ (\ref{iso-nu}) modulo $\h$ we conclude
that $b_0$ is invertible in $A$\,. Hence $b$ is
an invertible element of $\hA$\,.

On the other hand the cocycle condition for
$b$
$$
\nu^{-1}(a) * b - b * a =0\,, \qquad
a \in \hA
$$
implies that for every $a\in \hA$
$$
\nu^{-1}(a) = b * a * b^{-1}\,,
$$
where $b^{-1}$ is the inverse of $b$ in
the algebra $\hA$\,.

Thus, the automorphism $\nu$ is inner and
the implication $\Leftarrow$ follows. $\Box$

\section{Proof of Theorem \ref{ONA}.}
\label{Proof}
Let us pick on $X$ a nowhere vanishing
top degree form
$$
\om\in \G(\wedge^{\dim X} T^*X)
$$
and assign to this form the
modular vector field $v\in \h\G(TX)[[\h]]$
of the Poisson structure $\pi$
(\ref{pi-h1})\,. It is clear from
the definition (\ref{modular}) that
$v$ is a formal power series of
vector fields
\begin{equation}
\label{v}
v = \h v_1 + \h^2 v_2 + \dots
\end{equation}
starting with $\h v_1$ where $v_1$
is the modular vector field of the
initial Poisson structure $\pi_1$
on $X$\,.

Equation (\ref{pi-v}) implies that
$v$ is a cocycle in
the Poisson cochain complex of $\pi$\,.
Thus using the construction of
subsection \ref{q-deriv} we can lift
$v$ to a derivation
\begin{equation}
\label{cD-v}
\cD_v : \hA \to \hA
\end{equation}
of $\hA$ such that
\begin{equation}
\label{cD-v-h}
\cD_v = v \qquad m o d \quad \h^2\,.
\end{equation}

Let us consider the automorphism
of $\hA$ entering
the statement of Theorem \ref{ONA}
\begin{equation}
\label{phi}
\phi = \exp(\cD_v)
\end{equation}
and define an algebra structure
on $\hA[t,t^{-1}]$ by the following
rule:
\begin{equation}
\label{t}
(a t^m )\cdot (b t^k) =
a * \phi^m (b) t^{m + k}\,,
\end{equation}
where $a,b \in \hA$ and $m,k$ are
arbitrary integers.

Let us now analyze the bimodule structure on 
$$
\Ext^{\bul}_{(\hA [t,t^{-1}] )^{e}}
(\hA [t,t^{-1}], \hA [t,t^{-1}] \otimes \hA[t,t^{-1}] )
$$ 
using arguments of homological algebra.

For this we extend the action of the inverse 
\begin{equation}
\label{psi}
\psi = \phi^{-1}
\end{equation} 
of the automorphism
$\phi$ (\ref{phi}) to the Hochschild complex
$$
\Cbu(\hA, \hA\otimes \hA)
$$
in the natural way:
\begin{equation}
\label{psi-acts}
\psi(P)(a_1, \dots,  a_k) =
\phi^{-1}[ P(\phi(a_1), \dots, \phi(a_k)) ]
\end{equation}
where
$$
\phi^{-1}[b_1 \otimes b_2] = \phi^{-1}(b_1) \otimes \phi^{-1}(b_2)\,,
$$
and
$$
P \in C^k(\hA, \hA \otimes \hA)\,.
$$
The action (\ref{psi-acts}) is compatible
with the Hochschild coboundary operator and
with the internal $\hA$-bimodule structure, namely
\begin{equation}
\label{phi-bimod}
\psi(P) \cdot_{\intl} \phi^{-1}(a)  = \psi (P \cdot_{\intl} a)\,,
\qquad
\phi^{-1}(a) \cdot_{\intl} \psi(P)  = \psi (a \cdot_{\intl} P)\,.
\end{equation}

Thus, if we fix an isomorphism between
the $\hA$-bimodules
\begin{equation}
\label{isom-hA-nu}
\Ext^d(\hA, \hA\otimes \hA) \cong \hA \nu
\end{equation}
then the action
(\ref{psi-acts}) induces the following isomorphism
$\hpsi$
\begin{equation}
\label{hpsi}
\hpsi : \hA \nu \to \phi^{-1} \hA \nu \phi^{-1}
\end{equation}
from the $\hA$-bimodule $\hA \nu$ to
the $\hA$-bimodule\footnote{Hence $\phi \nu \phi^{-1}$
differs from $\nu$ by an inner automorphism.}
$\phi^{-1} \hA \nu \phi^{-1}$\,.

The isomorphism $\hpsi$ is uniquely determined by
the image of $1\in \hA \nu$\,. Let us denote this
image by $b_{\psi}$:
\begin{equation}
\label{b-psi}
b_{\psi} = \hpsi(1)\,.
\end{equation}
Since $\hpsi$ is an isomorphism the element
$b_{\psi}$ has to be invertible in $\hA$\,.

We claim that 
\begin{pred}
\label{esche-kak}
Given the isomorphism (\ref{isom-hA-nu})
we may construct the isomorphism
of $\hA[t,t^{-1}]$-bimodules
\begin{equation}
\label{ask}
\Ext^{\bul}_{(\hA [t,t^{-1}] )^{e}}
(\hA [t,t^{-1}], \hA [t,t^{-1}]\otimes \hA [t,t^{-1}] )
\cong
\begin{cases}
\hA [t,t^{-1}] \tnu  & {\rm if} \quad \bul = d + 1\,,\\
0  & {\rm otherwise}\,,
\end{cases}
\end{equation}
where
\begin{equation}
\label{tnu}
\tnu(a) = \nu(a)\,, \qquad \forall~~ a \in \hA \,,
\end{equation}
and
\begin{equation}
\label{tnu1}
\tnu(t) =  t \,  b_{\psi} \,.
\end{equation}
The isomorphism (\ref{ask}) is compatible 
with the natural action of the Euler 
field (\ref{Euler}).
\end{pred}
{\bf Proof.}
We remark that $\hA[t,t^{-1}]$ is
the twisted group algebra of the group
$\bbZ$\,. Thus, as a $\hA[t,t^{-1}]$-bimodule,
$\hA[t, t^{-1}]$ admits the following free
resolution:
\begin{equation}
\label{sPashei}
\Cbd(\bbZ, \Cbd(\hA, \hA[t,t^{-1}] \otimes \hA[t,t^{-1}]))\,,
\end{equation}
where $\Cbd(\bbZ, \,\,)$ denotes the standard
chain complex of the group $\bbZ$, the (right)
$\bbZ$-module structure on the chains of
the complex
$\Cbd(\hA, \hA[t,t^{-1}] \otimes \hA[t,t^{-1}])$
is given by
$$
(a t^n , a_1, \dots, a_k, b t^m) \, t  =
(a t^{n+1} , \phi^{-1}(a_1),  \dots,
\phi^{-1}(a_k), \phi^{-1}(b) t^{m-1})\,,
$$
$$
a, a_1, \dots, a_k, b \in \hA
$$
and for the definition of the Hochschild
chain complex $\Cbd(\hA, \hA[t,t^{-1}] \otimes \hA[t,t^{-1}])$
we use the internal $\hA$-bimodule structure in
$\hA[t,t^{-1}] \otimes \hA[t,t^{-1}]$\,.

The acyclicity of the resolution (\ref{sPashei})
is proved in the beginning of section $3$ in
\cite{Pasha} for an arbitrary discrete group
acting on an arbitrary associative algebra
with unit.

Using this resolution we conclude that
the $\hA$-bimodules
\begin{equation}
\label{oni}
\Ext^{\bul}_{(\hA [t,t^{-1}] )^{e}}
(\hA [t,t^{-1}], \hA [t,t^{-1}] \otimes \hA[t,t^{-1}] )
\end{equation}
can be computed as the total cohomology
of the following double complex
\begin{equation}
\label{double}
\Cbu(\bbZ, \Cbu(\hA, \hA[t,t^{-1}]
\otimes \hA[t,t^{-1}]))\,,
\end{equation}
where $\Cbu(\bbZ, \,\,)$ denotes the standard
cochain complex of the group $\bbZ$, the (right)
$\bbZ$-module structure on the cochains of
the complex
$\Cbu(\hA, \hA[t,t^{-1}] \otimes \hA[t,t^{-1}])$
is given by
\begin{equation}
\label{action-T}
(P \, T )(a_1, \dots, a_k)=
t \cdot_{\ext} P (\phi^{-1}(a_1), \dots,
\phi^{-1}(a_k))  \cdot_{\ext} t^{-1}\,, 
\end{equation}
$$
P\in C^k(\hA, \hA[t,t^{-1}] \otimes \hA[t,t^{-1}])\,,
$$
where $T$ denotes the (right) action of the generator 
of $\bbZ$\,, $\cdot_{\ext}$ denotes the multiplication
with respect to the external $\hA[t,t^{-1}]$-bimodule 
structure on $\hA[t,t^{-1}] \otimes \hA[t,t^{-1}]$\,, and the same 
external bimodule structure is used in the definition 
of the Hochschild cochain complex
$\Cbu(\hA, \hA[t,t^{-1}] \otimes \hA[t,t^{-1}])$\,.

Trivial $\bbZ$-module $\bbC$ admits the following 
simple resolution
$$
0 \to \bbC[t, t^{-1}]
\,\stackrel{\cdot (1-t)}{\longrightarrow}\, \bbC[t, t^{-1}] 
 \,\stackrel{t=1}{\longrightarrow}\,
\bbC \to 
 0\,.
$$
Therefore the cohomology of (\ref{double})
can be computed using the simpler double complex 
\begin{equation}
\label{s-double}
\begin{array}{ccccc}
\stackrel{\pa}{\to} & 
C^k(\hA, \hA[t,t^{-1}] \otimes \hA[t,t^{-1}]) &
\stackrel{\pa}{\to} & C^{k+1}(\hA, \hA[t,t^{-1}] \otimes \hA[t,t^{-1}]) &
\stackrel{\pa}{\to}  \\[0.3cm]
~ & \uparrow^{\, 1-T} & ~ & \uparrow^{\, 1-T} 
& ~  \\[0.3cm]
 \stackrel{\pa}{\to} & C^k(\hA, \hA[t,t^{-1}] \otimes \hA[t,t^{-1}]) &
\stackrel{\pa}{\to} & C^{k+1}(\hA, \hA[t,t^{-1}] \otimes \hA[t,t^{-1}]) &
\stackrel{\pa}{\to}  
\end{array}
\end{equation}
This double-complex is bounded from the left and 
it has length $2$ in the vertical direction.
 
Since the Hochschild complex
$\Cbu(\hA, \hA[t,t^{-1}] \otimes \hA[t,t^{-1}])$
splits into the following direct sum 
$$
\Cbu(\hA, \hA[t,t^{-1}] \otimes \hA[t,t^{-1}]) = 
\bigoplus_{n, m \in \bbZ} 
\Cbu(\hA, \hA t^n \otimes t^m \hA)
$$
Equation (\ref{k-neravno-d}) implies that the 
rows of the double-complex (\ref{s-double})
are acyclic beyond the dimension $d$\,. 
Using this observation it is not hard to show that 
$$
HH^{k} ( \hA[t,t^{-1}], \hA[t,t^{-1}] \otimes \hA[t,t^{-1}]) = 0
$$
if $k \neq d+1$ and
$$
HH^{d+1} ( \hA[t,t^{-1}], \hA[t,t^{-1}] \otimes \hA[t,t^{-1}]) 
\cong 
$$
\begin{equation}
\label{as-hA-mod}
 HH^{d} ( \hA, \hA[t,t^{-1}] \otimes \hA[t,t^{-1}])
\bigg/
[1-\hat{T}]\, \Big( HH^{d} ( \hA, \hA[t,t^{-1}] \otimes \hA[t,t^{-1}])
\Big)\,,
\end{equation}
where $\hat{T}$ is the action on $HH^{\bul}(\hA, \hA[t,t^{-1}] \otimes \hA[t,t^{-1}])$ induced by (\ref{action-T})\,. 

On the other hand, 
$HH^{d} ( \hA, \hA[t,t^{-1}] \otimes \hA[t,t^{-1}])$ is 
a free $\bbZ$-module generated by 
$HH^{d} ( \hA, \hA \otimes \hA[t,t^{-1}]) $ and hence
$$
HH^{d+1} ( \hA[t,t^{-1}], \hA[t,t^{-1}] \otimes \hA[t,t^{-1}]) 
\cong \hA[t,t^{-1}] \nu
$$
as an $\hA$-bimodule and a left module over 
$\hA[t,t^{-1}]$.

To determine the right module structure over 
$\hA[t,t^{-1}]$ we remark that for every cochain
$P\in C^k(\hA, \hA \otimes \hA)$ viewed as 
a cochain in  $C^k(\hA, \hA[t,t^{-1}] \otimes \hA[t, t^{-1}])$
\begin{equation}
\label{identity}
P \cdot_{\intl} t = ( t  \cdot_{\intl} \psi (P)) \, T \,.
\end{equation}
Due to (\ref{as-hA-mod}) Equation (\ref{identity})
implies that the right 
$\hA[t,t^{-1}]$-module structure of 
$$
HH^{d+1} ( \hA[t,t^{-1}], \hA[t,t^{-1}] \otimes \hA[t,t^{-1}])
$$ 
is indeed the one in (\ref{esche-kak}). 

It is obvious that the automorphism $\tnu$ (\ref{tnu}), 
(\ref{tnu1}) commutes with the action of the Euler 
field $\Eu$ (\ref{Euler})\,. Thus the field $\Eu$ acts 
of the $\hA[t,t^{-1}]$-bimodule $\hA[t,t^{-1}]\tnu$
in the natural way. It is clear from the 
construction that the isomorphism (\ref{ask}) is 
compatible with this action and this
completes the proof of Proposition \ref{esche-kak}. $\Box$

Due to Proposition \ref{dlia-hA-t} the DGLA module
$(\Cbu(\hA[t,t^{-1}]), \Cbd(\hA[t,t^{-1}]))$ of Hochschild 
(co)chains of the algebra $\hA[t,t^{-1}]$ with the product
(\ref{t}) is quasi-isomorphic to the DGLA module $(\Tp(X\times \bbC^{\times})[[\h]]_{restr}\,,\, \cAb(X\times \bbC^{\times})[[\h]]_{restr})$ where the DGLA  
$\Tp(X\times \bbC^{\times})[[\h]]_{restr}$ carries 
the differential $[\pi^t , \,\, ]$\,, the module 
$\cAb(X\times \bbC^{\times})[[\h]]_{restr})$ carries the 
differential $\cL_{\pi^t}$ and the Poisson bivector 
$\pi^t$ is given by the formula
\begin{equation}
\label{pi-t}
\pi^t  = \pi +  t\pa_{t} \wedge v\,.
\end{equation}

It is not hard to see that if $\om$ is the 
volume form we picked at the beginning of this 
section then the form
\begin{equation}
\label{om-t}
\om_{t} = \frac{d t}{t^2} \wedge \om
\end{equation}
satisfies the equation
\begin{equation}
\label{eq-om-t}
\cL_{\pi^t} \om_t =0\,,
\end{equation}
which implies that the Poisson structure 
(\ref{pi-t}) on $X\times \bbC^{\times}$ is 
unimodular.

Due to Proposition \ref{dlia-hA-t}, we may 
apply Theorem \ref{unimod-teo} to the algebra 
$\hA[t,t^{-1}]$ and deduce that
\begin{equation}
\label{vse-OK-0}
\Ext^{k}_{(\hA [t,t^{-1}] )^{e}}
(\hA [t,t^{-1}], \hA [t,t^{-1}] \otimes \hA[t,t^{-1}] )
\cong 0
\end{equation}
if $k \neq d+1$ and 
\begin{equation}
\label{vse-OK}
\Ext^{d+1}_{(\hA [t,t^{-1}] )^{e}}
(\hA [t,t^{-1}], \hA [t,t^{-1}] \otimes \hA[t,t^{-1}] )
\cong
 \hA [t,t^{-1}]
\end{equation}
as  $\hA [t,t^{-1}]$-bimodules.

Let us recapitulate the relevant part of the 
proof of Theorem \ref{unimod-teo} keeping track of 
the action of the Euler field $\Eu$ (\ref{Euler}).

Equation (\ref{eq-om-t}) implies that 
the form $\om_t$ (\ref{om-t}) is a cycle in the 
complex $\cAb(X\times \bbC^{\times})[[\h]]_{restr})$
with the differential $\cL_{\pi^t}$\,.
 
Due to Proposition \ref{dlia-hA-t} the homology 
of the complex $\cAb(X\times \bbC^{\times})[[\h]]_{restr})$
with the differential $\cL_{\pi^t}$ is isomorphic to the 
Hochschild homology $HH_{\bul}(\hA[t,t^{-1}], \hA[t,t^{-1}])$
of $\hA[t,t^{-1}]$\,. Furthermore, this isomorphism 
is compatible with the action of $\Eu$\,. Thus we may 
pull the class represented by $\om_t$ to a class in 
$$
c\in HH_{d+1}(\hA[t,t^{-1}], \hA[t,t^{-1}])\,.
$$

Since $\cL_{\Eu} \om_t = - \om_t$ we have
\begin{equation}
\label{Eu-c}
\Eu(c) = - c\,.
\end{equation}
Using Van den Bergh duality we 
pull $c$ to a class  
$$
c' \in HH^0(\hA[t,t^{-1}],  \tnu^{-1}\, \hA[t,t^{-1}])\,,
$$
where the automorphism $\tnu$ is defined in
(\ref{tnu}) and (\ref{tnu1}).

Due to Proposition \ref{esche-kak} Van den Bergh 
isomorphism 
$$
HH_{d+1}(\hA[t,t^{-1}],  \hA[t,t^{-1}]) \cong
HH^0(\hA[t,t^{-1}],  \tnu^{-1}\, \hA[t,t^{-1}])
$$
is compatible with the action of $\Eu$\,. Hence, 
\begin{equation}
\label{Eu-c1}
\Eu(c') = - c'\,.
\end{equation}

Thus if $b\in \hA[t,t^{-1}]$ represents the class 
$c'$ then
\begin{equation}
\label{Eu-b}
\Eu(b) = - b\,.
\end{equation}
In other words, the element $b$ is of the form  
\begin{equation}
\label{b}
b = b_0 t^{-1}\,, \qquad b_0\in \hA\,.
\end{equation}

The same argument as in the proof of Theorem 
\ref{unimod-teo} shows that $b_0$ is invertible.
While Equation (\ref{tnu}) and the cocycle 
condition for $b$
$$
\tnu^{-1}(a) * b - b * a =0\,, \qquad
\forall \quad a \in \hA
$$
imply that
$$
\nu^{-1}(a) = b_0 * \phi^{-1}(a) * b_0^{-1}\,, 
\qquad  \forall \quad a \in \hA\,,
$$
where the inverse $b_0^{-1}$ is taken 
in the algebra $\hA$\,.

Thus the automorphism $\nu$ differs from 
$\phi$ (\ref{phi}) by an inner automorphism and
Theorem \ref{ONA} is proved. $\Box$

\section{Concluding remarks}
Following V. Ginzburg \cite{Vitya} the
algebra $A$ of functions on the smooth affine
variety $X$ with trivial canonical bundle
gives us an example of a Calabi-Yau algebra.
(See Definition 3.2.3. and Corollary 3.3.2 in
\cite{Vitya}). Due to remark 3.2.8 in \cite{Vitya},
Theorem \ref{unimod-teo} implies that
the deformation quantization algebra $\hA$
is a Calabi-Yau algebra if and only if
a representative (\ref{pi-h1}) of
Kontsevich's class of $\hA$ is
unimodular\footnote{
It is not hard to see that the property
of being unimodular does not depend on
the choice of the representative.}.
This reformulation of Theorem \ref{unimod-teo}
is similar to the Kontsevich-Soibelman
conjecture about Calabi-Yau structures on an $A_{\infty}$
algebra $\mA$ (see conjecture 10.2.8 in \cite{K-Soi1}).
We should mention that structures very similar
to the Calabi-Yau structures \cite{K-Soi1}
were introduced and discussed in
\cite{DV}, \cite{DV1} for a large class
of connected graded algebras.

It should be remarked that the results
of G. Felder and B. Shoikhet \cite{F-Sh}
also show a special role of unimodular
Poisson manifolds in deformation quantization.
More precisely in \cite{F-Sh} the authors
prove a statement closely related to
the cyclic formality conjecture \cite{TT}, 
\cite{Tsygan}. Using this
result they show that Kontsevich's star-product with
the harmonic angle function is cyclic. Finally,
they also prove a generalization of the
Connes-Flato-Sternheimer
conjecture \cite{CFS} on closed star-products
in the Poisson case.
It would be interesting to find a relationship
between the results of \cite{F-Sh} and Theorem
\ref{unimod-teo} proved in this paper.

We would like to mention recent preprint
\cite{quadrat} by S. Launois and L. Richard
in which the authors consider the algebra
$B_{\h}$ of functions on
a uniparametrized
quantum affine space and notice that, in this case,
the modular class of the corresponding Poisson
bracket can be restored as
the semiclassical limit of the modular
automorphism $\nu$ (\ref{Ext-d}) of this
algebra $B_{\h}$\,. Theorem \ref{ONA} shows
that, in general, the semiclassical limit
of the modular automorphism may not restore
the modular class. In particular, it is
easy to construct an example in which the
semiclassical limit of the modular automorphism vanishes,
while the modular class is still non-zero.

Theorem \ref{unimod-teo} shows that, in the
unimodular case, we have the Van den Bergh
duality \cite{VB} isomorphism:
\begin{equation}
\label{dual-HH}
V : HH^{\bul}(\hA, \hA)\,
\stackrel{\sim}{\to}\, HH_{d -\bul}(\hA, \hA)\,,
\end{equation}
where $d$ is the dimension of $X$\,.

On the other hand, due to \cite{K} we have
an isomorphism between
Hochschild cohomology $HH^{\bul}(\hA, \hA)$
of $\hA$ and the Poisson cohomology (\ref{HP})
of $\pi$ (\ref{pi-h1}). Furthermore, due
to \cite{thesis} and \cite{Sh} we also have
an isomorphism between
Hochschild homology $HH_{\bul}(\hA, \hA)$
of $\hA$ and the Poisson homology (\ref{HPd})
of $\pi$ (\ref{pi-h1}). This observation poses
a question about a relationship between
the Van den Bergh duality (\ref{dual-HH})
and the duality of P. Xu \cite{Xu}
(see Theorem 4.8 in \cite{Xu}).

In the general
(non-unimodular) case a version of Xu's duality
(see, for example, Eq. (10) in \cite{quadrat})
involves the Poisson homology with
coefficients in a module over the Poisson
algebra. A generalization of this
duality statement to Lie algebroids
was proposed in \cite{ELW} and extended further
to the framework of module categories
in \cite{Block}.
 It would be interesting
to investigate a relation of this duality
to Van den Bergh duality for the corresponding
deformation quantization algebra.

We should also mention
paper \cite{BZh}, in which K.A. Brown and J.J. Zhang
discuss the Van den Bergh duality for
noetherian Hopf algebras. They showed that,
in this case, the dualizing module is
also determined by an outer automorphism of
the algebra. They call this outer automorphism
{\it the Nakayama automorphism}.

A natural generalization of the modular class
to the case $Q$-manifolds \cite{Tomsk} shows up
in the quantization of gauge systems. More precisely,
it was shown in \cite{Tomsk} that this class happens
to be the first obstruction to the existence
of a quantum master action in the BV quantization.

Finally, we would like to mention papers
\cite{DF} and \cite{Yvetta} in which the authors
consider the hierarchy generated by the
modular vector field of a Poisson-Nijenhuis
manifold and relate it to a
hierarchy of bihamiltonian vector fields.
It would be interesting to investigate
quantum versions of these hierarchies
and their relation to quantization
of bihamiltonian systems \cite{Olver},
\cite{MAO}.

\section*{Appendix. Poisson, Hamiltonian and
log-Hamiltonian vector fields.}
In this section we introduce the notion
of log-Hamiltonian 
vector fields and prove
some of their useful properties. In paper \cite{BW}
these vector fields are called ``integral derivations''. 
However, we believe that the name ``log-Hamiltonian'' 
is more appropriate for such vector fields. 
Very similar objects show up in literature
in various contexts \cite{AAM}, \cite{DK}, \cite{Goto}, 
\cite{OR}.

Let $X$ be a smooth affine variety and
$\pi$ be the Poisson structure
(\ref{pi-h}). A vector field
$ w \in \G(TX)[[\h]]$ is called a
{\it Poisson vector field} of $\pi$ if
\begin{equation}
\label{Poisson}
[\pi, w]_{SN} = 0\,.
\end{equation}
It is called a {\it Hamiltonian vector field}
if there exists a function
$f\in A[[\h]]$ such that
\begin{equation}
\label{Hamil}
w = [\pi, f]_{SN}\,.
\end{equation}
Finally, we say that $w$ is a {\it log-Hamiltonian
vector field } of $\pi$ if there exists
a unit $f$ of $A[[\h]]$ such that
\begin{equation}
\label{log-Ham}
w = f^{-1} [\pi, f]_{SN}\,.
\end{equation}
Notice that, the units of the ring
$A[[\h]]$ are the formal power series
in $\h$
$$
f = f_0 + \h f_1 + \h^2 f_2 + \dots\,,
\qquad f_i \in A
$$
starting with a nowhere vanishing function
$f_0\in A$\,.

The Poisson
vector fields are degree $0$ cocycles
and the Hamiltonian vector fields
are degree $0$ coboundaries
in the complex (\ref{Poisson-com})\,.
Every log-Hamiltonian vector
field is a Poisson vector field.

Let us list some simple properties
of the log-Hamiltonian vector fields
in the following proposition:
\begin{pred}
\label{proper-log}
~~~
\begin{enumerate}
\item A linear combination of
log-Hamiltonian vector fields with
integer coefficients is again a
log-Hamiltonian vector field.

\item If $w_1$ is a Poisson vector
field and $w_2$ is a log-Hamiltonian
vector field then their Lie bracket
$[w_1, w_2]$ is a Hamiltonian vector
field.

\item If $g\in \h A[[\h]]$ then the Hamiltonian
vector field $w = [\pi, g]_{SN}$
is a log-Hamiltonian vector field.

\end{enumerate}
\end{pred}
~\\
{\bf Proof.}
\begin{enumerate}
\item follows from the equations
$$
n f^{-1}[\pi,f]_{SN} = f^{-n}[\pi, f^n]_{SN}\,,
\qquad f_1^{-1}[\pi,f_1]_{SN} +
 f_2^{-1}[\pi,f_2]_{SN} =  (f_1 f_2)^{-1}
 [\pi, f_1 f_2]_{SN}\,,
$$
where $n$ is an integer, $f, f_1,$ and $f_2$ are
units of $A[[\h]]$\,.

\item Using the fact that $w_1$ is
a Poisson vector field it is not hard
to prove that
$$
[w_1, f^{-1}[\pi, f]_{SN} ]_{SN} =
[\pi, f^{-1} w_1(f)]_{SN}\,,
$$
which immediately implies $2$\,.

\item follows from the obvious equation
$$
e^{-g}[\pi, e^g]_{SN} = [\pi, g]_{SN}\,,
$$
where $e^{\pm g}$ makes sense because
$g\in \h A[[\h]]$\,. $\Box$
\end{enumerate}

This proposition shows that the log-Hamiltonian
vector fields form a lattice in the space
of Poisson vector fields. The classes of
the log-Hamiltonian vector fields in
$HP^1(X, \pi)$ form a lattice in the center
of the Lie algebra $HP^1(X, \pi)$\,.

In subsection \ref{q-deriv} we construct
a map $\cD$ (\ref{cD})
from the $\bbC[[\h]]$-module of Poisson
vector fields to the $\bbC[[\h]]$-module of
the derivations of the corresponding
deformation quantization algebra $\hA$\,.

If a Poisson vector field $w$ starts with
$\h$ then so does the corresponding
derivation $\cD_{w}$ and it makes sense
to speak about the automorphism $\exp(\cD_{w})$\,.
Let us prove that
\begin{teo}
\label{lift-log-Ham}
If
\begin{equation}
\label{w-f}
w = - f^{-1}[\pi, f]_{SN}
\end{equation}
for some unit $f\in A[[\h]]$ then
there is a function $\tf$ such that
\begin{equation}
\label{f-tf}
\tf = f \qquad  mod \quad \h\,,
\end{equation}
and
\begin{equation}
\label{w-tf}
\exp(\cD_{w})(a) = \tf * a * \tf^{-1}\,,
\qquad \forall \quad a\in \hA\,,
\end{equation}
where $\tf^{-1}$ is the inverse of $\tf$
in the algebra $\hA$\,.
\end{teo}
{\bf Proof.} In the proof we use the construction
and notation from subsection \ref{q-deriv}.
It is also useful
to understand the construction of the
maps $\la_T$ and $\la_D$ from \cite{thesis}
(see Propositions 13 and 14 in chapter 4 of
\cite{thesis}).

It is obvious that the lift
$\la_T(f)$ satisfies the following property
\begin{equation}
\label{la-f}
\la_T(f) = f (1+g)\,,
\end{equation}
where $g$ is a section of $\SX$ satisfying
the equation
\begin{equation}
\label{g-y-0}
g|_{y^i=0} = 0\,.
\end{equation}
Therefore, it makes sense to speak about
the section
\begin{equation}
\label{tg}
\tg = \ln(1+g)\,,
\end{equation}
defined via the Taylor
power series of $\ln(1+x)$ around
the point $x=0$\,.

Equation (\ref{w-f}) implies that
\begin{equation}
\label{la-w-f}
\la_T(w) = - \la_T(f)^{-1}[\la_T(\pi), \la_T(f)]_{SN}\,.
\end{equation}

Furthermore, since $\la_T(f) = f (1+g)$
and $f$ does not depend on the fiber
coordinate $y$
$$
\la_T(w) = - (1+g)^{-1}[\la_T(\pi), (1+g)]_{SN}\,,
$$
or equivalently, in terms of $\tg$ (\ref{tg})
\begin{equation}
\label{la-w-tg}
\la_T(w) = - [\la_T(\pi), \tg \,]_{SN}\,.
\end{equation}

Since $\la_T(f)$ is a flat
section with respect to the Fedosov
connection (\ref{DDD}), $D \la_T(f) =0$,
we have
\begin{equation}
\label{D-tg}
D \tg = - f^{-1} d f\,,
\end{equation}
where $d$ is the de Rham differential.

Combining (\ref{la-w-tg}) and (\ref{D-tg})
we get
\begin{equation}
\label{D-tg1}
D \tg + [\la_T(\pi), \tg]_{SN}
= - \la_T(w) - f^{-1} d f\,,
\end{equation}
Hence,
\begin{equation}
\label{cK-tg}
D \tcK_1(\tg) +
\pa_{\dia} \tcK_1(\tg) =
- \tcK_1(\la_T(w)) - \tcK_1(f^{-1} d f)\,.
\end{equation}
Since $f^{-1}d f$ is independent of
the fiber coordinates $y^i$ we have
$$
\tcK_1(f^{-1} d f) = f^{-1} d f\,.
$$
Therefore,
\begin{equation}
\label{cK-tg1}
D \tcK_1(\tg) +
\pa_{\dia} \tcK_1(\tg) =
- \tcK_1(\la_T(w)) - f^{-1} d f \,.
\end{equation}
In subsection \ref{q-deriv}
the cochain $\tcK_1(\la_T(w))$ was
denoted by $W$ and the derivation $\cD_{w}$
was defined by Equation (\ref{cD-def}).
Using this equation we rewrite (\ref{cK-tg1})
as
\begin{equation}
\label{cK-tg11}
D \tcK_1(\tg) +
\pa_{\dia} \tcK_1(\tg) =
- \la_D(\cD_w) - (D + \pa_{\dia}) \Phi(W_1)
- f^{-1} d f\,,
\end{equation}
where $W_1$ is the component of $W$ of exterior
degree $1$\,.

Since $\tg$ is a cochain of
$(\OmT[[\h]], [\la_T(\pi),\,])$
of degree $-1$\,, 
$$
\tcK_1(\tg)\in \Om^0(\SX)[[\h]]\,.
$$
Hence, combing the components of exterior
degree $0$ and $1$ in (\ref{cK-tg11}), we get
\begin{equation}
\label{cD-w}
\la_D(\cD_w) = - \pa_{\dia} (\tcK_1(\tg) + \Phi(W_1))
\end{equation}
\begin{equation}
\label{D-G}
D (\tcK_1(\tg) + \Phi(W_1)) = - f^{-1} d f
\end{equation}
Since the series $\pi$ (\ref{pi-h})
belongs to $\h \, T^1_{poly}(X)[[\h]]$
and $\tg\Big|_{y^i=0}= 0$ due to (\ref{g-y-0})
$$
\tcK_1(\tg) \Big |_{\h= y^i =0} = 0\,.
$$
Furthermore, using the equation defining
the homotopy operator $\Phi$ (see Eq. (4.36)
in chapter 4 of \cite{thesis}) we conclude that
$$
\Phi(W_1)\Big|_{y^i} = 0\,.
$$
Thus the element
\begin{equation}
\label{G}
G = \tcK_1(\tg) + \Phi(W_1) \in \G(\SX)[[\h]]
\end{equation}
satisfies the property
\begin{equation}
\label{G-y-0}
G\Big|_{\h=y^i=0} = 0\,,
\end{equation}
and hence, the $\dia$-exponent of
$G$
\begin{equation}
\label{exp-G}
\exp_{\dia}(G) = 1+
\sum_{k\ge 1} \frac{1}{k!}
\underbrace{G\dia G \dia \dots \dia G}_{k\,{\rm times}}
\end{equation}
makes sense.

Equation (\ref{cD-w}) implies that
$$
\exp (\la_D(\cD_w))(b) =
\exp_{\dia}(G)\dia b \dia \exp_{\dia}(-G)\,,
\qquad \forall\quad b\in \G(\SX)[[\h]]
$$
and, since $f$ does not depend on fiber
coordinates $y^i$,
\begin{equation}
\label{exp}
\exp (\la_D(\cD_w))(b) =
f \, \exp_{\dia}(G)\dia b \dia f^{-1}\, \exp_{\dia}(-G)\,,
\qquad \forall\quad b\in \G(\SX)[[\h]]\,.
\end{equation}
Let us prove that
\begin{equation}
\label{D-f-G}
D ( f \, \exp_{\dia}(G)) = 0\,.
\end{equation}

Indeed,
$$
D ( f \, \exp_{\dia}(G)) = d f \, \exp_{\dia}(G)
+ f  \sum_{k\ge 1} \frac1{k!}
\sum_{p+q = k-1}
\underbrace{G\dia G \dia \dots \dia G}_{p\,{\rm times}}
\dia (DG) \dia
\underbrace{G\dia G \dia \dots \dia G}_{q\,{\rm times}}\,.
$$
Equation (\ref{D-G}) implies that $D G = - f^{-1} d f $.
Therefore,
$$
D ( f \, \exp_{\dia}(G)) = d f \, \exp_{\dia}(G)
- f \sum_{k\ge 1} \frac1{k!}
\sum_{p+q = k-1}
\underbrace{G\dia G \dia \dots \dia G}_{p\,{\rm times}}
\dia (f^{-1} d f) \dia
\underbrace{G\dia G \dia \dots \dia G}_{q\,{\rm times}}\,.
$$
But $f^{-1} d f$ does not depend on the tangent
coordinates $y^i$, and hence belongs to the center
of the algebra $(\G(\SX)[[\h]], \dia)$\,.
Thus, the latter equation gives (\ref{D-f-G}).

Equation (\ref{D-f-G}) implies that
the element $f \, \exp_{\dia}(G)$
of $\Om^0(C^{-1}(\SX))[[\h]]$ belongs to
the image of the embedding $\la_D$
from (\ref{chain})\,. In other words
there exists $\tf\in A[[\h]]$ such that
\begin{equation}
\label{tf1}
\la_D(\tf) = f \, \exp_{\dia}(G)\,.
\end{equation}

In fact $\tf$ can be obtained by
simply setting $y^i=0$ in
$f \, \exp_{\dia}(G)$
\begin{equation}
\label{tf}
\tf = f \, \exp_{\dia}(G) \Big|_{y^i = 0}
\end{equation}

Using (\ref{diamond}),
(\ref{exp}) and (\ref{tf1}) we immediately
conclude that $\tf$ satisfies equation
(\ref{w-tf}).
Equation (\ref{f-tf}) follows from the fact
that the element $G$ (\ref{G}) satisfies
(\ref{G-y-0})\,. This completes the proof of Theorem
\ref{lift-log-Ham}. $\Box$

Using claims $2$ and $3$ of Proposition
\ref{proper-log}, the Baker-Campbell-Hausdorff
formula, and Theorem \ref{lift-log-Ham}
it is not hard to prove that
\begin{cor}
\label{sense}
If $w\in \h\G(TX)[[\h]]$
is Poisson vector field and
$w_{\log}$ is a log-Hamiltonian vector
field then the automorphism
$$
\exp(\cD_{w}+ \cD_{w_{\log}})\circ \exp(-\cD_{w})
$$
of $\hA$ is inner. $\Box$
\end{cor}

}

~\\

\noindent\textsc{Department of Mathematics,
Northwestern University,
Evanston, IL 60208 \\
\emph{E-mail address:} {\bf vald@math.northwestern.edu}}

\end{document}